\documentclass[12pt]{article}
\usepackage[numbers,sort&compress]{natbib}
\usepackage{amssymb,amsmath,mathrsfs,amsfonts,amsthm}
\usepackage{enumerate}
\usepackage{paralist}
\usepackage{color}
\usepackage{epsfig}
\usepackage{graphicx}
\usepackage{setspace}
\usepackage{appendix}
\begin{document}
 \def\pd#1#2{\frac{\partial#1}{\partial#2}}
\def\dfrac{\displaystyle\frac}
\let\oldsection\section
\renewcommand\section{\setcounter{equation}{0}\oldsection}
\renewcommand\thesection{\arabic{section}}
\renewcommand\theequation{\thesection.\arabic{equation}}

\newtheorem{thm}{Theorem}[section]
\newtheorem{cor}[thm]{Corollary}
\newtheorem{lem}[thm]{Lemma}
\newtheorem{prop}[thm]{Proposition}
\newtheorem*{con}{Conjucture}
\newtheorem*{questionA}{Question}
\newtheorem*{thmA}{Theorem A}
\newtheorem*{thmB}{Theorem B}
\newtheorem{remark}{Remark}[section]
\newtheorem{definition}{Definition}[section]

\title{Nonlocal to Local Convergence of Stefan Problems Under Optimal Convergence Condition
\thanks{M. Zhou was partially supported by  NSF of China (No. 12271437).  F. Li  was partially supported by  NSF  of China (No. 12371213, 12126609).} }

\author{
Xinfu Chen    {\thanks{Email:  xinfu@pitt.edu}}\\
{Department of Mathematics, University of Pittsburgh,}\\
{\small Pittsburgh, PA 15260, USA.}\\
Fang Li {\thanks{ E-mail: lifang55@mail.sysu.edu.cn}}
\\
{School of Mathematics, Sun Yat-sen University,}\\
{\small No. 135, Xingang Xi Road, Guangzhou 510275, P. R. China.} \\
Maolin Zhou {\thanks{Corresponding author. Email:  zhouml123@nankai.edu.cn}}\\
{Chern Institute of Mathematics and LPMC, Nankai University,}\\
{\small Tianjin 300071, P. R. China.} }

\date{}
\maketitle{}

\begin{abstract}
In this paper, we consider a  free boundary problem with a nonlocal diffusion kernel function $k(x)$. Due to the long distance exchange effect of nonlocal diffusion, the free boundary can expand discontinuously, which makes the problem rather complicated. Among other things, we propose the optimal convergence condition without  assuming the symmetry or compactness of  $k$, i.e., the Fourier transform of $k$  satisfies
$$\hat{k}(\xi)=1-|\xi|^2+o(|\xi|^2)\ \ \mbox{ as }\xi\rightarrow 0,$$
and discover  an equivalent characterization of this optimal condition. More importantly, by the employment of the variational inequality, the apriori estimates and the Fourier transform, we demonstrate that, along a series of properly rescaled kernel functions, the corresponding solutions to the nonlocal free boundary problems converge to the solution of the classical Stefan problem under the proposed optimal condition.
\end{abstract}



\section{Introduction}

The {\it classical Stefan problem} is well known to describe the evolution of the interface between two phases of a substance undergoing a phase change, for example the melting of a solid, such as ice to water. {\it Latent heat}, defined as the heat or energy that is absorbed or released during a phase change of a substance, acts as an energy source or sink at a moving solid-liquid interface, and the resulting boundary condition is known as the {\it Stefan boundary condition}.

In this paper, we  propose and study   {\it the nonlocal version of one-phase Stefan problem}
\begin{equation}\label{one-nonlocalstefan}
\begin{cases}
\displaystyle\gamma_t(t,x)= d \int_{\{\gamma>0\}} \hspace{-0.2cm} k(x-y)\gamma(t,y)dy- d \gamma(t,x)\chi_{\{\gamma>0\}}& t>0, \ x\in\mathbb R^n,\\
\gamma(0,x)=\gamma_0(x)  & x\in \mathbb{R}^n,
\end{cases}
\end{equation}
where  $\chi_E$ denotes the characteristic function of $E$, the kernel function $k$ satisfies
\begin{description}
\item[(K)]$\ \ \ \ \ \ \ \ \ \ \ \ \ \ \ \ \ \ \ \ \ \ \ \ \displaystyle k\in C(\mathbb R^n)\cap L^\infty(\mathbb R^n),\; k\geq 0,\; k(0)>0,~\int_{\mathbb{R}^n}k(x)dx=1.$
\end{description}
This nonlocal diffusion operator is widely used to model  a long range diffusion process and  appears commonly in different types of models in ecology.
See \cite{Allen1996, HMMV, Kot1996, Lee2001,  Medlock2003,  Meysman2003, Mogilner-E} and the references therein.

In addition to ecology, nonlocal operators are also widely used to  model many other applied situations involving nonlocal interactions, and they have rich variations and extensions.
For a systematic introduction of modeling and analysis, see the books \cite{Andreu-J-Rossi-T} by Andreu-Vaillo, Maz\'{o}n, Rossi, and Toledo-Melero and \cite{Du-book} by Du and the survey paper \cite{DuICM} by Du. In the development of rigorous theoretical framework for nonlocal models, the connection between local models and their nonlocal analogue is one of the fundamental questions.
In particular, the convergence relations between different types of evolution equations with local and nonlocal diffusions have been thoroughly elaborated in the book \cite{Andreu-J-Rossi-T}.
{\it Among other things, the main purpose of this paper is to demonstrate the nonlocal to local convergence of one-phase Stefan problem under the optimal condition.}

\medskip
For clarity, we always assume that
\begin{equation}\label{initialdata}
\Omega_0 \ \textrm{is a smooth and bounded domain in}\ \mathbb R^n,\    \ell_0 \ \textrm{is a positive constant,}
\end{equation}
and for the initial data,
\begin{equation}\label{initial-onephase}
\gamma_0(x)\in L^{\infty}(\mathbb R^n),\ \gamma_0(x)= -\ell_0  \ \ \textrm{for} \  x\in \mathbb{R}^n \setminus \bar\Omega_0,\ \gamma_0|_{\bar\Omega_0}\geq 0, \ \gamma_0|_{\bar\Omega_0} \not\equiv 0.
\end{equation}
Also denote
$$
\gamma^+(t,x)= \gamma(t,x)\chi_{\{\gamma>0\}}.
$$
This will be used whenever it is more convenient.

To better elaborate  the formulation of the model, we  first consider the classical {\it one-phase Stefan problem}, which is the description, typically, of the melting of a body of ice, maintained at zero degree centigrade, in contact with a region of water initially in $\Omega_0$. Based on latent heat and conservation of energy, the model is formulated as follows
\begin{equation}\label{localstephan-n}
\begin{cases}
\displaystyle \theta_t(t,x)= d\Delta\theta (t,x)  & t>0,\ x\in \{\theta(t, \cdot)>0\}, \\
\nabla_x \theta \cdot \nabla_x s= -\ell_0   &t>0,\ x\in\partial \{\theta(t, \cdot)>0\},  \\
\theta=0 &t>0,\ x\in\partial \{\theta(t, \cdot)>0\},  \\
\theta(0,x)=u_0(x)   & x\in \bar\Omega_0,
\end{cases}
\end{equation}
where $\theta =\theta (t,x)$ denotes the water's temperature, the free boundary $\partial \{\theta(t, \cdot)>0\}$ at the time $t$ is given by the equation $s(x)=t$. Also  set $s(x)=0$ if $x\in \bar\Omega_0$. There are many famous papers on the regularity of the free boundary, such as \cite{C,CPS,F}.

On the basis of latent heat,  the nonlocal version of one-phase  Stefan problem (\ref{one-nonlocalstefan}) is proposed  and the essence is at the time $t$,
\begin{itemize}
\item if $x \in \{ x\in\mathbb R^n \, | \,  \gamma(t, x) \leq 0 \}$, then  only it can absorb energy from outside;
\item if $x\in  \{  x\in\mathbb R^n \, | \,  \gamma(t, x) > 0 \}$, then not only  it can absorb energy from outside, but also transfer its energy outside.
\end{itemize}
Here, the value $\ell_0$ is in the status of latent heat, $\gamma$ equal to $-\ell_0$ corresponds to the status of ice at  zero degree centigrade and $\gamma$ reaching zero  represents that there has already been sufficient energy accumulated here for the phase change.

The nonlocal version of  one-phase Stefan problem was  also proposed and studied in \cite{BCQuira2012}.  Some discussions  will be placed when the results obtained in this paper are related to those derived in \cite{BCQuira2012}.
Moreover, the fractional  Stefan problem was treated in \cite{ACaffarelli2010}, and more general, the  Stefan problem with anomalous diffusion was investigated in \cite{ACaffarelliM2022}.

Notice that, in \cite{BCQuira2012}, besides the assumption \textbf{(K)}, it is always assumed that the kernel function is  compactly supported and radially symmetric. However, in this paper, we carry on the investigation of    the  nonlocal  one-phase Stefan  problem  (\ref{one-nonlocalstefan}) under the  assumption \textbf{(K)}.

\medskip


First of all,  we establish  results about the local existence and global existence  for the nonlocal Stefan problems.

\begin{thm}\label{thm-wellposedness}
Assume that in the problem  (\ref{one-nonlocalstefan}), the kernel function $k$  satisfies  the assumption  \textbf{(K)},  the condition (\ref{initialdata}) is valid and the initial data satisfies (\ref{initial-onephase}).    Then the problem (\ref{one-nonlocalstefan})   admits a unique classical solution $\gamma(t,\cdot) \in L^{\infty}(\mathbb R^n)$ defined for all $t>0$, and  $\gamma$ satisfies the estimate
\begin{equation}\label{thm-wellposedness-bdd}
\mathrm{ess} \inf_{\mathbb R^n} \gamma_0  \leq \gamma(t,x)  \leq  \mathrm{ess} \sup_{\mathbb R^n} \gamma_0   \ \ \ \textrm{for} \ t>0,\, x\in\mathbb R^n.
\end{equation}
Moreover,  if $\gamma_0|_{\bar\Omega_0}\in C(\bar\Omega_0)$, then $\gamma(t, \cdot)$ is continuous in $\bar\Omega_0$ and $\mathbb{R}^n \setminus \bar\Omega_0$ for any $t>0$.
\end{thm}

Next, we  demonstrate some  fundamental properties related to {\it expansion, boundedness and continuity} of  the free boundaries  $\partial\{ x\in\mathbb R^n \ |\ \gamma(t, x)\geq 0 \}$ in the nonlocal one-phase  Stefan problem (\ref{one-nonlocalstefan}).
In particular, we prove that  if the initial domain $\Omega_0$ is convex and the kernel function $k(x)$ is radially symmetry and strictly decreasing on $r=|x|$, then the free boundary will continuously expand; otherwise, {\it the jumping phenomenon}, i.e., the discontinuous expansion of the free boundary may happen due to certain combination of the initial data and the kernel function. Some typical examples are constructed for elaboration.

We emphasize that the jumping phenomenon
never appears in the classical Stefan problem.
In fact, since the nonlocal diffusion describes the  movement between non-adjacent spatial locations, the jumping phenomena is naturally expected and reflects the essential differences between local and nonlocal diffusion operators. This reveals not only the complexity of  evolution of free boundaries in nonlocal Stefan problems but also the difficulties in verifying the nonlocal to local convergence of Stefan problems.

Now, we discuss the convergence relations between local and nonlocal  Stefan problems.
For simplicity, for $\epsilon>0$, denote
$$
k_\epsilon(x)={1\over \epsilon^n} k({x \over\epsilon}),\ \eta_\epsilon(x)= {1\over \epsilon^n}\eta({x \over\epsilon}).
$$
Before we present the main results, we briefly explain what should be the natural and optimal assumptions on the nonlocal kernel functions in the studies of convergence relations  between models with local and nonlocal diffusions.
Define the Fourier transform of the kernel function  $k$ as follows
$$
\displaystyle\hat k (\xi) = \int_{\mathbb R^n} e^{-ix\cdot \xi} k(x)dx.
$$
Based on the properties of the Fourier transform, one observes that   for $\phi\in L^1(\mathbb R^n) \bigcap C^2(\mathbb R^n)$
$$
\int_{\mathbb R^n} e^{-ix\cdot \xi}  \left({1\over \epsilon^2}\int_{\mathbb R^n} k_{\epsilon}(x-y)\phi(y)dy- {1\over \epsilon^2}\phi(x)\right) dx
=
{1\over \epsilon^2}\left( \hat k(\epsilon \xi) -1 \right)\hat \phi(\xi),
$$
$$
\int_{\mathbb R^n} e^{-ix\cdot \xi} \Delta\phi(x) dx = -|\xi|^2 \hat \phi(\xi),
$$
and for fixed $\xi$,
$$
\lim_{\epsilon\rightarrow 0}{1\over \epsilon^2}\left( \hat k(\epsilon \xi) -1 \right)\hat \phi(\xi)= -A|\xi|^2 \hat \phi(\xi)
$$
under the condition
\begin{equation}\label{assumption-hat k}
\hat k (\xi )  = 1- A|\xi |^2 +o(|\xi |^2) \ \ \ \textrm{as}  \  \xi\rightarrow 0,
\end{equation}
where $A>0$ is a constant.
This observation  indicates that {\it the condition (\ref{assumption-hat k}) should be  optimal  in the studies of nonlocal approximation of Laplacian operator.}



In the following,  we  discover  an important equivalent characterization of the optimal convergence  condition (\ref{assumption-hat k}), where neither the radial symmetry nor the compact supportedness property is required.

\begin{prop}\label{prop-kernel}
Assume that $k$ satisfies the assumption  \textbf{(K)}.
Then the following two statements are equivalent.
\begin{itemize}
\item[(i)] For $1\leq j, h\leq n,\  j\neq h$, $\displaystyle \int_{\mathbb R^n} x_jk(x)dx =0,\  \int_{\mathbb R^n} x_jx_h k(x)dx =0, \ \int_{\mathbb R^n} x_j^2k(x)dx ={1\over n}\int_{\mathbb R^n} |x|^2 k(x) dx<+\infty.$
\item[(ii)] The Fourier transform of $k$  satisfies the assumption (\ref{assumption-hat k}), where
$$
A={1\over 2n}\int_{\mathbb R^n} |x|^2 k(x) dx.
$$
\end{itemize}
\end{prop}

The equivalent characterization of the optimal convergence  condition
(\ref{assumption-hat k}) obtained in Proposition \ref{prop-kernel}(i) clearly describes the requirements on the decay rate of the kernel function $k$ when $|x|$ is large, and its asymmetry.
In addition, if the kernel function $k$ is radially symmetric,  Proposition \ref{prop-kernel}(i) is reduced to
$$
\int_{\mathbb R^n} |x|^2 k(x) dx<+\infty,
$$
and  obviously, radially symmetric kernel functions with compact support satisfy Proposition \ref{prop-kernel}(i). We also point out that the  decay rate of the kernel function is an important indicator in studying the qualitative properties of  models with nonlocal diffusions. For example, it is well known that in the nonlocal analogue of the Fisher-KPP equation, the decay rate of the kernel function at infinity determines the existence of traveling wave solutions. See  \cite{Yagisita2009}  and the references therein.

In the book \cite{Andreu-J-Rossi-T}, the nonlocal to local convergence of Cauchy problems is verified for radially symmetric kernel functions satisfying the expected optimal  condition \eqref{assumption-hat k}.
However, for both   Dirichlet and Neumman  boundary problems, the convergence relations are proved under the assumptions  that the kernel functions are radially symmetric and compactly supported. See \cite{Rossi-Dirichlet, Rossi-Neumann} and the references therein.  It remains unclear what should be the optimal convergence conditions for problems with Dirichlet or Neumann boundary condition.

In this paper,  the convergence relations between  local and  nonlocal one-phase Stefan problems  are verified under   the optimal condition (\ref{assumption-hat k}) without additional assumptions about symmetry or compactness of kernel functions.

To begin with, we rescale the problem (\ref{one-nonlocalstefan})  as follows
\begin{equation}\label{one-nonlocalstefan-epsilon}
\begin{cases}
\displaystyle \gamma_{\epsilon t}(t,x)=  {1\over \epsilon^2}\int_{\mathbb R^n} k_{\epsilon}(x-y)\gamma_{\epsilon}^+(t,y)dy- {1\over \epsilon^2}\gamma_{\epsilon}^+(t,x)   &t>0,\ x\in\mathbb{R}^n, \\
\gamma_{\epsilon}(0,x)=\gamma_0(x)    & x\in \mathbb{R}^n,
\end{cases}
\end{equation}
and denote
$$
\gamma_{\epsilon}^+(t,x)= \gamma_{\epsilon}(t,x)\chi_{\{\gamma_{\epsilon}(t,x)>0\}}.
$$

\begin{thm}\label{thm-convergence}
In the problem (\ref{one-nonlocalstefan-epsilon}), assume that the kernel function satisfies  the assumption  \textbf{(K)},  the condition (\ref{initialdata}) is valid and the initial data satisfies (\ref{initial-onephase}). Also, assume that  the Fourier transform of $k$ satisfies   the condition (\ref{assumption-hat k}).
Then for any given $T>0$, $\gamma_{\epsilon}^+$ converges to the solution $\theta$ of the one-phase Stefan problem (\ref{localstephan-n}) as $\epsilon \rightarrow 0$ in the following sense:
$$
 \int_0^t \gamma_{\epsilon}^+(\tau,x) d\tau \rightarrow \int_{\min \{s(x),t\}}^t \theta(\tau, x)d\tau \ \ \textrm{a.e. in } (0,T) \times \mathbb R^n,
$$
where we set $d=A$ and $u_0=\gamma_0|_{\bar\Omega_0}$ in the problem (\ref{localstephan-n}).
\end{thm}


In general, in the process of convergence of the nonlocal Stefan problem (\ref{one-nonlocalstefan-epsilon}) as $\epsilon \rightarrow 0$ , the  main difficulties arising from the following aspects:
    \begin{itemize}
\item The domain where $\gamma_\epsilon (t,x)$ is nonnegative is always changing as the time $t$ and the parameter $\epsilon $ vary.
\item It is well known that  problems with nonlocal diffusion operators lack  spatial regularity.
\item It is possible that the free boundary of the nonlocal Stefan problem evolves discontinuously, i.e., the jumping phenomena happens.
\item It is delicate to estimate the amount of  energy accumulation based on the decay rate of the kernel function $k$ when $|x|$ is large.
\end{itemize}

The nonlocal to local convergence of one-phase Stefan problems is also studied in \cite[Theorem 5.2]{BCQuira2012}.  Since the proof is based on Taylor expansion, the authors confirm the convergence relations between weak solutions  under  the additional conditions that the kernel function is radially symmetric and  compactly supported.

Theorem \ref{thm-convergence} verifies  the convergence relations between classical solutions under the optimal condition (\ref{assumption-hat k}).
To achieve this improvement, we  borrow the idea of the variational inequality for  the classical Stefan problem to  transform the nonlocal Stefan  problem into a whole space problem, and then combine the apriori estimates derived in Section 2.3 and  the Fourier transform to derive certain convergence of solutions.

This paper is organized as follows.  Theorem \ref{thm-wellposedness} and some useful apriori estimates for the  nonlocal problem (\ref{one-nonlocalstefan}) are established in Section 2. In Section 3, we focus on the proof of Proposition \ref{prop-kernel}, which demonstrates the equivalent characterization of the optimal convergence condition (\ref{assumption-hat k}). Section 4 is devoted to  the proof of Theorem  \ref{thm-convergence}, which is about  the convergence relations between local and nonlocal Stefan problems.  In Appendixes A and B, we prove some fundamental properties related to  expansion, boundedness and continuity of  free boundaries.


%

\section{Wellposedness and preliminaries}

\subsection{Local and global existence}

We first verify the local and global existence to the the nonlocal version of the one-phase Stefan problem (\ref{one-nonlocalstefan}).
\begin{proof}[Proof of Theorem \ref{thm-wellposedness}]
Denote
$M_0 = \| \gamma_0\|_{L^{\infty}(\mathbb R^n)} $,  $\mathbb Y=  L^{\infty}(\mathbb R^n), $
 for $s>0$,
$$
\mathbb X_s =\left\{      \phi \in C([0,s), \mathbb Y)\, \big| \, \phi(0,\cdot)=\gamma_0(\cdot),\  \|\phi(t,\cdot)\|_{L^{\infty}(\mathbb R^n)}\leq 2M_0,  \, t\in [0,s) \right\},
$$
and
$$
\| \phi\|_{C([0,s),\mathbb Y)}= \sup_{0\leq t<s} \|\phi(t,\cdot)\|_{L^{\infty}(\mathbb R^n)}.
$$
For $\phi \in \mathbb X_s$, $0<t<s$, define
\begin{equation*}
\begin{aligned}
\mathcal{T} \phi &= \gamma_0(x) + d \int_0^t\int_{\{\phi >0\}} k(x-y)\phi(\tau,y)dyd\tau- d \int_0^t \phi(\tau,x)\chi_{\{\phi>0\}}d\tau.
\end{aligned}
\end{equation*}
Then it is routine to show that  $\mathcal{T} \phi\in C([0,s),\mathbb Y)$, $\mathcal{T} \phi (0,\cdot) = \gamma_0(\cdot)$ and
\begin{eqnarray*}
\|\mathcal{T} \phi \|_{C([0,s),L^{\infty}(\mathbb R^n))} \leq M_0 + 2d s\| \phi\|_{C([0,s),\mathbb Y)}\leq  M_0 + 4 d s M_0.
\end{eqnarray*}
Moreover, for $\phi_1, \phi_2\in \mathbb X_s$,
\begin{eqnarray*}
\|\mathcal{T} \phi_1 - \mathcal{T} \phi_2 \|_{C([0,s),\mathbb Y)} \leq 2 d s\| \phi_1 -  \phi_2 \|_{C([0,s),\mathbb Y)} .
\end{eqnarray*}
Thus it is obvious that there exists $t_0>0$, which depends $a,\, b$ and $M_0$ only and is sufficiently small, such that for $0<s \leq  t_0$, $\mathcal{T}$ maps $ \mathbb X_s$ into $ \mathbb X_s$ and $\mathcal{T}$ is a contraction mapping in $ \mathbb X_s$. Hence by  the contraction mapping theorem, 
for $0<s \leq  t_0$, there exists a unique $\gamma \in \mathbb X_s$  satisfying
\begin{equation*}
\begin{aligned}
\gamma(t,x) &= \gamma_0(x) + d\int_0^t\int_{\{\gamma >0\}} k(x-y)\gamma(\tau,y)dyd\tau- d\int_0^t \gamma(\tau,x)\chi_{\{\gamma>0\}}d\tau
\end{aligned}
\end{equation*}
for $0<t<s$, $x\in\mathbb R^n$. Thus, obviously  $\gamma$  is the unique  solution to the problem (\ref{one-nonlocalstefan}).

Let $(0,T_{\max})$ denote the maximal time interval for which the solution $\gamma(t,x)$  of the problem (\ref{one-nonlocalstefan}) exists.
It remains to show  $T_{\max}= + \infty$. For this purpose,  it suffices to show that  $\| \gamma(t,\cdot) \|_{L^{\infty}(\mathbb R^n)}$ is bounded in $(0,T_{\max})$.
To be more specific, we claim that {\it $\gamma$ satisfies  the estimate (\ref{thm-wellposedness-bdd}) in $(0,T_{\max})$.}

Fix any $0<T<T_{\max}$.
First, assume that the kernel function $k$ is compactly supported.   Then since  $\bar\Omega_0$ is bounded, it is standard to show that $\{\gamma(t, x)\geq 0\}$ remains bounded for $0<t<T$.
Based on the problem (\ref{one-nonlocalstefan}), for any $1<p< +\infty$, $0<t<T$,  one has
$$
(\gamma^+)^{p-1} \gamma^+_t(t,x)\leq d (\gamma^+)^{p-1}\left(  \int_{\{\gamma>0\}} \hspace{-0.2cm} k(x-y)\gamma(t,y)dy- \gamma(t,x)\chi_{\{\gamma>0\}}\right).
$$
Then direct computation yields that for $0<t<T$,
\begin{eqnarray}
&& {1\over p} {d\over dt}\int_{\mathbb R^n}(\gamma^+(t,x))^pdx \cr
&\leq&  d \int_{\mathbb R^n}(\gamma^+(t,x))^{p-1}\left(  \int_{\mathbb R^n}  k(x-y)\gamma^+(t,y)dy- \gamma^+(t,x) \right)dx\cr
&\leq &  d  \int_{\mathbb R^n}(\gamma^+(t,x))^{p-1}  \left(\int_{\mathbb R^n} k(x,y)dy\right)^{p-1\over p}\left(\int_{\mathbb R^n}k(x,y)(\gamma^+(t,y))^pdy\right)^{1\over p} dx- d \|\gamma^+(t,\cdot)\|_{L^p(\mathbb R^n)}^p\cr
&\leq & d \|\gamma^+(t,\cdot)\|_{L^p(\mathbb R^n)}^{p-1}\left(\int_{\mathbb R^n}\int_{\mathbb R^n}k(x,y)(\gamma^+(t,y))^pdydx\right)^{1\over p}- d\|\gamma^+(t,\cdot)\|_{L^p(\mathbb R^n)}^p\leq 0.\nonumber
\end{eqnarray}
Hence  for any $1<p< +\infty$, $0<t<T$,
$$
\|\gamma^+(t,\cdot)\|_{L^p(\mathbb R^n)}\leq \|\gamma^+(0,\cdot)\|_{L^p(\mathbb R^n)},
$$
and it follows that
$$
\gamma(t,x)  \leq  \mathrm{ess} \sup_{\mathbb R^n} \gamma_0,  \ \   0<t<T,\, x\in\mathbb R^n.
$$
The lower bound on $\gamma$ that
$$
\gamma(t,x)  \geq  \mathrm{ess} \inf_{\mathbb R^n} \gamma_0, \ \   0<t<T,\, x\in\mathbb R^n
$$
is obvious.
The claim is verified for compactly supported kernel functions since $T\in(0, T_{\max})$ is arbitrary.

Now we consider the case that the kernel function $k$ is not compactly supported.  Then there exists a series of kernel functions $k_j$,  $j\geq 1$, which are compactly supported, satisfy the assumption \textbf{(K)}, and
\begin{equation}\label{kernel-L1}
\lim_{j\rightarrow \infty}\| k_j -k\|_{L^1(\mathbb R^n)}=0.
\end{equation}
Let $\gamma_j$ denote the solution to the problem (\ref{one-nonlocalstefan}) with $k$ replaced by $k_j$. Set $w_j= \gamma -\gamma_j$, $j\geq 1$. Then $w_j$ satisfies
\begin{equation*}
\begin{cases}
\displaystyle (w_j)_t(t,x) =d \int_{\{\gamma>0\}} k(x-y)\gamma(t,y)dy- d\gamma(t,x)\chi_{\{\gamma>0\}}\\
\displaystyle \hspace{2.3cm}  -d \int_{\{\gamma_j>0\}} k_j(x-y)\gamma_j(t,y)dy+ d\gamma_j(t,x)\chi_{\{\gamma_j>0\}},& t>0, \ x\in\mathbb R^n, \\
w_j(0,x)= 0, & x\in \mathbb{R}^n.
\end{cases}
\end{equation*}
Then for  $w_j> 0$, direct computation yields that
\begin{eqnarray}
(w_j)_t(t,x) &\leq& d\int_{\{\gamma>0\}} k(x-y) \left(\gamma(t,y)-\gamma_j(t,y)\right)dy  + d \int_{\{\gamma>0\}} k(x-y) \gamma_j(t,y) dy\cr
&& -d \int_{\{\gamma_j>0\}}\left( k_j(x-y)-k(x-y)\right)\gamma_j(t,y)dy- d \int_{ \{\gamma_j>0\}} k(x-y) \gamma_j(t,y) dy\cr
&\leq& d \| w_j\|_{L^\infty(\mathbb R^n)} + d M_0 \| k_j -k\|_{L^1(\mathbb R^n)},\nonumber
\end{eqnarray}
where the last inequality follows from    the fact that $\gamma_j$ satisfies the estimate  (\ref{thm-wellposedness-bdd}).
Similarly for $w_j<0$, we have
$$
(-w_j)_t(t,x) \leq  d\| w_j\|_{L^\infty(\mathbb R^n)} + dM_0 \| k_j -k\|_{L^1(\mathbb R^n)}.
$$
The above two inequalities indicate that for $0<t<T_{\max}$,
\begin{eqnarray}\label{pf-thm-wj}
&& | w_j(t,x) | = \lim_{\delta \rightarrow 0^+} \int_0^t {\partial\over \partial\tau }\left[w_j^2(\tau, x)+\delta^2\right]^{1\over 2}\, d\tau = \lim_{\delta \rightarrow 0^+} \int_0^t \frac{w_j(\tau, x)}{\left[w_j^2(\tau, x)+\delta^2\right]^{1\over 2}}{\partial\over \partial\tau } w_j(\tau, x) d\tau\cr
&\leq& d \int_0^t  \| w_j\|_{L^\infty(\mathbb R^n)}(\tau) d\tau + d M_0 \| k_j -k\|_{L^1(\mathbb R^n)}t.
\end{eqnarray}
Denote
$$
h_j(t) = \int_0^t  \| w_j\|_{L^\infty(\mathbb R^n)}(\tau) d\tau,
$$
then (\ref{pf-thm-wj}) implies that for $0<t<T_{\max}$,
$$
h_j'(t) \leq d h_j(t)+ d M_0 \| k_j -k\|_{L^1(\mathbb R^n)}t.
$$
Direct computation yields that for $0<t<T_{\max}$,
$$
h_j(t) \leq {1\over d}M_0 e^{dt}\| k_j -k\|_{L^1(\mathbb R^n)},
$$
which, together with (\ref{pf-thm-wj}), indicates that for $0<t<T_{\max}$,
\begin{equation}\label{pf-thm-ex1}
 \| w_j\|_{L^\infty(\mathbb R^n)}(t)\leq M_0\left(e^{dt}+ d t\right)\| k_j -k\|_{L^1(\mathbb R^n)}.
\end{equation}
Therefore,  (\ref{kernel-L1}) and the fact that $\gamma_j$ satisfies  the estimate (\ref{thm-wellposedness-bdd}) for all $j\geq 1$ imply the desired claim for general kernel functions under the assumption \textbf{(K)}.

At the end, it is routine to verify that  if $\gamma_0|_{\bar\Omega_0}\in C(\bar\Omega_0)$, then $\gamma(t, \cdot)$ is continuous in $\bar\Omega_0$ and $\mathbb{R}^n \setminus \bar\Omega_0$ for any $t>0$.
\end{proof}

\subsection{The maximum principle}
We first present the comparison principle for the nonlocal version of  one-phase Stefan problem (\ref{one-nonlocalstefan}) and omit the proof since it is standard.

\begin{prop}\label{prop-cp}
Assume that the conditions of Theorem \ref{thm-wellposedness} are valid. Also assume that $\gamma^*_0\in L^{\infty}(\mathbb R^n)$,  $\gamma^*_0(x)= -\ell_0$  for  $x\in \mathbb{R}^n \setminus \bar\Omega_0$.  Let $\gamma^*$ denote the solution to the problem (\ref{one-nonlocalstefan}) with initial data $\gamma^*_0$.  If  $\gamma^*_0 \geq \gamma_0$, then $\gamma^* \geq\gamma$ for all $t>0$.
\end{prop}

Moreover, we establish a type of strong maximum principle for the nonlocal version of one-phase Stefan problem (\ref{one-nonlocalstefan}).

\begin{prop}\label{prop-strongmp}
Under the conditions of Theorem \ref{thm-wellposedness}, for given $s\geq 0$, we have $\gamma(t,x)>0$ in $\Omega(s)$ when $t > s$, where
$$
\Omega(s) = \{ x\in\mathbb R^n \ |\ \gamma(s, x)\geq 0 \}.
$$
\end{prop}

\begin{proof}
First, we claim that {\it if $x\in\{x\in\Omega(s) \, | \, \gamma(s,x)>0\}$, then for $t>s$, $\gamma(t,x) >0 $.}
Due to the continuity of the solution in  $t$, we only need consider the case that
$s>0$.
According to
$$
\gamma_t(t,x)= d \int_{\{\gamma>0\}}  k(x-y)\gamma(t,y)dy- d \gamma(t,x)\chi_{\{\gamma>0\}}\geq - d \gamma(t,x)\chi_{\{\gamma>0\}},
$$
the claim follows immediately.

Next we consider the initial domain $\bar\Omega_0$.   Set
$$
\gamma_{0\delta}(x)=
\begin{cases}
\gamma_0(x)+\delta  & x\in\bar\Omega,\\
\gamma_0(x)     & x\in \mathbb R^n\setminus\bar\Omega_0,
\end{cases}
$$
where $\delta>0$, and let $\gamma_{\delta}$ denote the solution to the problem (\ref{one-nonlocalstefan}) with the initial data (\ref{initial-onephase}), where $\gamma_0$ is replaced by $\gamma_{0\delta}$.  Thanks to the above claim, one sees that  $\gamma_{\delta}(t,x)>0$ for $t>0$, $x\in \bar\Omega_0$. By letting $\delta\rightarrow 0^+$, it is routine to derive that $\gamma(t, x)\geq 0$ for $t>0$, $x\in\bar\Omega_0$, i.e., $\bar\Omega_0\subseteq \Omega(t)$ for $t>0$.

Moreover, since $\gamma_0|_{\bar\Omega_0}\geq 0, \ \gamma_0|_{\bar\Omega_0} \not\equiv 0$, the  claim  at the beginning  indicates that
$\{ x\in\bar\Omega_0 \,| \, \gamma(t, x)>0\}$ is not empty for $t>0$.
Suppose that there exists $t_0>0$ such that $\gamma(t_0, x)$ touches zero somewhere in $\bar\Omega_0$.  By choosing
$$
x_0\in \partial  \{ x\in\bar\Omega_0 \,| \, \gamma(t_0, x)>0\} \bigcap\{ x\in\bar\Omega_0 \,| \, \gamma(t_0, x)=0\},
$$
we have
$$
0\geq \gamma_t(t_0,x_0)= d \int_{\{\gamma(t_0, y)>0\}} k(x_0-y)\gamma(t_0,y)dy >0,
$$
where the strict inequality is due to the assumption \textbf{(K)} and the choice of $x_0$. This is a contradiction and thus $\gamma(t, x)>0$ for $t>0$, $x\in\bar\Omega_0$.

It remains to consider the  set $\{x\in\Omega(s) \setminus\bar\Omega_0 \, | \, \gamma(s,x)=0\}$, when it is not empty.   Fix $x^*\in\{x\in\Omega(s) \setminus\bar\Omega_0 \, | \, \gamma(s,x)=0\}$
and let $s_1$ denote the moment when $\gamma(t,x^*)$ first touches zero. Obviously $s_1\leq s$ and by the equation satisfied by $\gamma$, we have
$$
\ell_0 = d \int_0^{s_1} \int_{\{\gamma(t,y)>0\}}  k(x^*-y)\gamma(t,y)dydt.
$$
Then obviously there exists $t_1\in(0,s_1)$ such that
\begin{equation}\label{pf-prop-positive}
\int_{\{\gamma(t_1,y)>0\}}  k(x^*-y)\gamma(t_1,y)dy>0.
\end{equation}
We claim that {\it for any $t>t_1$,  $\displaystyle\int_{\{\gamma(t,y)>0\}}  k(x^*-y)\gamma(t,y)dy>0$.}
Suppose that the claim is not true, i.e., there exists $t_2>t_1$ such that
$$
\int_{\{\gamma(t_2,y)>0\}}  k(x^*-y)\gamma(t_2,y)dy = 0.
$$
This implies that $\gamma(t_2,y)\leq 0$ in the set $\{y\in\mathbb R^n \, | \, k(x^*-y)>0\}$. Again thanks to the claim at the beginning, we have $\gamma(t_1,y)\leq 0$ in the set $\{y\in\mathbb R^n \, | \, k(x^*-y)>0\}$, which contradicts to (\ref{pf-prop-positive}). The claim is proved.

According to this claim and the choice of $s$, $x^*$, one sees that
$$
\gamma_t(s,x^*)= d \int_{\{\gamma(s,y)>0\}}  k(x^*-y)\gamma(s,y)dy>0.
$$
Hence  for $t>s\geq 0$, $\gamma(t,x^*)>0$.
\end{proof}

\subsection{Apriori estimates}
In this subsection, we derive some apriori estimates  in the nonlocal Stefan problem (\ref{one-nonlocalstefan}), which    are useful in proving convergence relations between  local and  nonlocal  Stefan problems.

\begin{prop}\label{prop-Lp-estimate}
Under the assumptions of Theorem \ref{thm-wellposedness},  there exists  a constant $C_1>0$, which depends on the initial data only, such that for given $1\leq p\leq \infty$, we have
    $$
    \| \gamma^+(t, \cdot)\|_{L^p(\mathbb R^n)}\leq C_1 \mbox{ for }  t>0.
    $$
\end{prop}

\begin{proof}
Notice that if $\phi\in L^1(\mathbb R^n)\bigcap L^{\infty}(\mathbb R^n)$, then  for any $p>1$, $\phi \in L^p(\mathbb R^n)$ and
$$
\|\phi\|_{L^p(\mathbb R^n)}\leq \left(\|\phi\|^{p-1}_{L^{\infty}(\mathbb R^n)}\|\phi\|_{L^1(\mathbb R^n)} \right)^{1\over p}\leq \left(\|\phi\|_{L^{\infty}(\mathbb R^n)}+1\right)\left(\|\phi\|_{L^1(\mathbb R^n)}+1 \right).
$$
Hence it suffices to verify the statements   for $p=1$ and $p=\infty$.

Indeed, when $p=\infty$, the conclusion is obvious due  to Theorem \ref{thm-wellposedness}, i.e.,
\begin{equation}\label{pf-lem-u-infty-bdd}
\| \gamma(t, \cdot)\|_{L^{\infty}(\mathbb R^n)}\leq   \| \gamma_0\|_{L^{\infty}(\mathbb R^n)} .
\end{equation}

In order to estimate $\| \gamma^+(t, \cdot)\|_{L^1(\mathbb R^n)}$,
we first consider the case that $k$ is compactly supported. Let $\hat\gamma(t,x)$ denote the solution to the problem (\ref{one-nonlocalstefan}) with the initial data replaced by
$$
\hat\gamma(0,x)= \| \gamma_0|_{\bar\Omega_0}\|_{L^{\infty}(\Omega_0)},\ x\in\bar\Omega_0,\ \ \hat\gamma(0,x)=-\ell_0,\ x\in \mathbb R^n\setminus \bar\Omega_0.
$$
By Theorem \ref{thm-wellposedness} and Proposition \ref{prop-cp}, we have
\begin{equation}\label{pf-lem-hat}
\hat\gamma(t,x)\geq \gamma(t,x),\  -\ell_0\leq \hat\gamma(t,x)\leq \| \gamma_0|_{\bar\Omega_0}\|_{L^{\infty}(\Omega_0)},\  \ t>0, \,x\in\mathbb R^n.
\end{equation}
Since  $\bar\Omega_0$ is bounded and $k$ is compactly supported,  for $t>0$, it is routine to show that $\{\hat\gamma(t, x)\geq 0\}$   remains bounded.
Set
$$
\Sigma^+(t)= \bigcup_{0<\tau<t}\{\hat\gamma(\tau, x)\geq 0\}. 
$$ 
Then by direct computation, for $0<\tau<t$,
\begin{eqnarray*}
\int_{\Sigma^+ (t)} \hat\gamma_{ \tau}(\tau,x) dx
\leq d \int_{ \Sigma^+(t)} \int_{\mathbb R^n} k(x-y)  \hat\gamma^+ (\tau,y)dy dx - d \int_{ \Sigma^+ (t)}  \hat\gamma^+ (\tau,x) dx \leq 0.
\end{eqnarray*}
Thus
\begin{eqnarray*}
0\leq \int_{ \Sigma^+ (t)}\hat \gamma (t,x) dx  \leq \int_{ \Sigma^+ (t)} \hat\gamma(0,x) dx  = -\ell_0\ |\Sigma^+(t) \setminus \bar\Omega_0| + \| \gamma_0|_{\bar\Omega_0}\|_{L^{\infty}(\Omega_0)} \ |\bar \Omega_0|,
\end{eqnarray*}
which implies that
\begin{equation*}
|\{\hat\gamma (t, x)\geq 0\} |  \leq |\Sigma^+ (t) | \leq \left(1+ { \| \gamma_0|_{\bar\Omega_0}\|_{L^{\infty}(\Omega_0)}\over \ell_0} \right) |\bar \Omega_0|.
\end{equation*}
Hence  thanks to (\ref{pf-lem-hat}), for any given  $t>0$
\begin{equation}\label{pf-lem-domain-bdd}
\| \gamma^+(t, \cdot)\|_{L^1(\mathbb R^n)}\leq {\| \gamma_0|_{\bar\Omega_0}\|_{L^{\infty}(\Omega_0)}}\left(1+ { \| \gamma_0|_{\bar\Omega_0}\|_{L^{\infty}(\Omega_0)}\over \ell_0} \right) |\bar \Omega_0|  .
\end{equation}

Now consider the case that the kernel function $k $ satisfies the assumption \textbf{(K)}, but is not compactly supported.  Then there exists a series of kernel functions $k_j$, $j\geq 1$, which are compactly supported, satisfy the assumption \textbf{(K)}, and
$$
\lim_{j\rightarrow \infty}\| k_j -k_{\epsilon}\|_{L^1(\mathbb R^n)}=0.
$$
Let $\gamma_j$ denotes the solution to the problem (\ref{one-nonlocalstefan}) with $k$ replaced by $k_j$. Similar to the proof of (\ref{pf-thm-wj}) in the proof of Theorem \ref{thm-wellposedness}, we have
$$
\lim_{j\rightarrow \infty} \| \gamma_j^+ -\gamma^+ \|_{L^{\infty}(\mathbb R^n)}\leq \lim_{j\rightarrow \infty} \| \gamma_j -\gamma \|_{L^{\infty}(\mathbb R^n)}=0.
$$
This, together with  (\ref{pf-lem-domain-bdd}), implies that for any $R>0$,
$$
\int_{B_R(0)}\gamma^+ (t,x) dx = \lim_{j\rightarrow \infty} \int_{B_R(0)}\gamma_j^+ (t,x) dx \leq {\| \gamma_0|_{\bar\Omega_0}\|_{L^{\infty}(\Omega_0)}}\left(1+ { \| \gamma_0|_{\bar\Omega_0}\|_{L^{\infty}(\Omega_0)}\over \ell_0} \right) |\bar \Omega_0|.
$$
Since $R$ is arbitrary, for any given  $t>0$,
$$
\| \gamma^+(t, \cdot)\|_{L^1(\mathbb R^n)}  \leq {\| \gamma_0|_{\bar\Omega_0}\|_{L^{\infty}(\Omega_0)}}\left(1+ { \| \gamma_0|_{\bar\Omega_0}\|_{L^{\infty}(\Omega_0)}\over \ell_0} \right) |\bar \Omega_0|.
$$
The proof  is complete.
\end{proof}

\begin{prop}\label{prop-uniform-difference}
Under the assumptions of Theorem \ref{thm-wellposedness},    we have
    $$
\int_{\mathbb R^n}\left|\gamma(t,x+h) - \gamma(t,x)\right|dx \leq  \int_{\mathbb R^n}\left|\gamma_0(x+h) -  \gamma_0(x)\right|dx, \  \ t>0, \,x\in\mathbb R^n.
$$
\end{prop}


\begin{proof}
First of all, fix $x,\, h\in\mathbb R^n$.
For $\delta \neq 0$, introduce
$$
\mu_{\delta}(X) = \left(X^2+\delta^2\right)^{1\over 2}.
$$
According to the problem (\ref{one-nonlocalstefan}) satisfied by $\gamma$,  it is routine  to verify  that
\begin{eqnarray}
&& {\partial\over \partial t}\mu_{\delta}(\gamma(t,x+h) -  \gamma(t,x)) \cr
&=&\frac{\gamma(t,x+h) -  \gamma(t,x)}{\left[\left(\gamma(t,x+h) -  \gamma(t,x)\right)^2+\delta^2\right]^{{1\over 2}} } \left(\gamma(t,x+h) -  \gamma(t,x)\right)_t \cr
&=& \frac{\gamma(t,x +  h)- \gamma(t,x)}{\left[\left(\gamma(t,x + h) - \gamma(t,x)\right)^2+\delta^2\right]^{{1\over 2}} } \cr
&&\times  \left( d \int_{\mathbb R^n} [k(x-  y)  \gamma^+(t,y + h)- \gamma^+(t,y)] dy - d [ \gamma^+(t,x+ h)- \gamma^+(t,x)] \right) \cr
&\leq &\frac{|\gamma(t,x +  h)- \gamma(t,x)|}{\left[\left(\gamma(t,x + h) - \gamma(t,x)\right)^2+\delta^2\right]^{{1\over 2}} } \cr
&&\times  \left(d \int_{\mathbb R^n} k(x-  y)| \gamma^+(t,y + h)- \gamma^+(t,y) |dy - d | \gamma^+(t,x+ h)- \gamma^+(t,x) |\right),
\end{eqnarray}
which yields that
\begin{eqnarray*}
&& \left|\gamma (t,x+h) -  \gamma(t,x)\right| -\left|\gamma_0(x+h) -  \gamma_0(x)\right|\\
&=&\lim_{\delta  \rightarrow 0}\left[\mu_{\delta}(\gamma (t,x+h) -  \gamma(t,x))- \mu_{\delta}(\gamma_0(x+h) -  \gamma_0(x))\right]\\
&=& \lim_{\delta  \rightarrow 0} \int_0^t  {\partial\over \partial \tau}\mu_{\delta}(\gamma(\tau,x+h) -  \gamma(\tau,x)) d\tau\\
&\leq &  \int_0^t\left(d\int_{\mathbb R^n} k(x-y)| \gamma^+(\tau,y+h)-  \gamma^+(\tau,y) |dy - d | \gamma^+(\tau,x+h)- \gamma^+(\tau,x) | \right)d\tau.
\end{eqnarray*}
Thus for any $R>0$,
\begin{eqnarray}
&& \int_{B_R(0)}\left|\gamma (t,x+h) -  \gamma(t,x)\right|dx -\int_{B_R(0)}\left|\gamma_0(x+h) -  \gamma_0(x)\right|dx\cr
&\leq &  \int_0^t\int_{B_R(0)}\left(d\int_{\mathbb R^n} k(x-y)| \gamma^+(\tau,y+h)-  \gamma^+(\tau,y) |dy - d | \gamma^+(\tau,x+h)- \gamma^+(\tau,x) | \right) dx d\tau\cr
&\leq &  \int_0^t\left(d\int_{\mathbb R^n} | \gamma^+(\tau,x+h)-  \gamma^+(\tau,x) |dx -  d\int_{B_R(0)}| \gamma^+(\tau,x+h)-  \gamma^+(\tau,x) |dx \right)  d\tau.
\end{eqnarray}
Notice that $\gamma_0(x+h) -  \gamma_0(x)$ is compactly supported. Then thanks to Proposition \ref{prop-Lp-estimate}, by letting $R\rightarrow \infty$, we have for $t>0$,
$$
\int_{\mathbb R^n}\left|\gamma (t,x+h) - \gamma (t,x)\right|dx \leq  \int_{\mathbb R^n}\left|\gamma_0(x+h) -  \gamma_0(x)\right|dx.
$$
The proof is complete.
\end{proof}


\section{Fundamental properties of nonlocal Stefan problem}
In this section, we  investigate  some  fundamental properties related to {\it expansion, boundedness and continuity} of  free boundaries  $\partial\{ x\in\mathbb R^n \ |\ \gamma(t, x)\geq 0 \}$ in   the nonlocal version of one-phase Stefan problem (\ref{one-nonlocalstefan})
\begin{equation*}
\begin{cases}
\displaystyle \gamma_t(t,x)= d\int_{\mathbb R^n} k(x-y)\gamma^+(t,y)dy- d\gamma^+(t,x)   & t>0,\  x\in\mathbb{R}^n,\\
\gamma(0,x)=\gamma_0    & x\in \mathbb{R}^n.
\end{cases}
\end{equation*}

\subsection{Expansion, boundedness and continuity}
It is known that for the models with nonlocal diffusion, the regularity in space variables can not be improved. Hence, to study the evolution of the free boundaries, we always assume that the initial data $\gamma_0$ satisfies the extra condition that $\gamma_0|_{\bar\Omega_0}\in C(\bar\Omega_0)$.
For simplicity,  denote
\begin{equation}\label{omega}
\Omega(t) = \{ x\in\mathbb R^n \ |\ \gamma(t, x)\geq 0 \}
\end{equation}
where $\gamma(t,x)$ is the solution to the problem (\ref{one-nonlocalstefan}).

\begin{thm}[Expansion and boundedness]\label{thm-expansion-boundedness}
In the problem (\ref{one-nonlocalstefan}), assume that the conditions of Theorem \ref{thm-wellposedness} are valid  and the initial data $\gamma_0$ satisfies the extra condition that $\gamma_0|_{\bar\Omega_0}\in C(\bar\Omega_0)$.
We have the following statements.
\begin{itemize}
\item[(i)] \underline{\rm Expansion:} there exists $t_0>0$ such that $\Omega(t) = \Omega(0)$ for $0\leq t\leq t_0$ and $\Omega(t_1)\subseteq \Omega(t_2)$ for $0<t_1<t_2$.
\item[(ii)] \underline{\rm Boundedness:}  there exists $R>0 $, which depends on the initial data only,  such that $\Omega(t) \subseteq B_R(0) \ \textrm{for all} \ t>0$.
\end{itemize}
\end{thm}

We remark that Theorem \ref{thm-expansion-boundedness}(i)  is also proved in \cite{BCQuira2012}, where the kernel function is assumed to be compactly supported and radially symmetric.
The proof of Theorem \ref{thm-expansion-boundedness}(ii) is quite lengthy. In order not to interrupt the main theme of this paper, we move the proof of Theorem \ref{thm-expansion-boundedness} to Appendix A.

We further investigate the continuous expansion of $\Omega(t)$.
For convenience, we prepare an extra assumption about the kernel function as follows
 \begin{description}
\item[(K1)] $k(x)$ is  radially symmetric, decreasing in $|x|$.
\end{description}

The following result confirms the continuous expansion of $\Omega(t)$ under the extra conditions that the initial domain $\bar\Omega_0$ is convex and the kernel function $k$ satisfies \textbf{(K1)}.

\begin{thm}[Continuity]\label{thm-continuity}
Under the conditions of Theorem \ref{thm-wellposedness}, in the problem (\ref{one-nonlocalstefan}), if  we additionally  assume that the initial domain $\bar\Omega_0$ is convex and the  assumption \textbf{(K1)} is valid, then $\Omega(t)$ expands continuously.
\end{thm}

The proof of Theorem \ref{thm-continuity} will be included in Appendix B.

\subsection{Jumping phenomena}

The purpose of this subsection is to investigate a natural question arising from  Theorem \ref{thm-continuity}:

\smallskip

\noindent\textbf{Question}: What happens to the expansion of $\Omega(t)$ without these extra conditions that the initial domain $\bar\Omega_0$ is convex?
\medskip

We construct two examples to show that the population range could generate at a distant place when either the kernel function $k$ does not satisfy \textbf{(K1)} or the initial domain $\Omega_0$ is not convex.  This is so called {\it jumping phenomena}.

We also point out that,  if allowing the initial data to be  nonconstant outside $\bar\Omega_0$, similar to  \cite[Theorem 4.6]{BCQuira2012}, where the kernel function is assumed to be compactly supported and radially symmetric, jumping phenomena could happen by choosing initial data properly. Indeed, the conclusion is also valid as long as the kernel function satisfies the assumption \textbf{(K)}. We omit the proof since it is similar. In the following two examples, we keep the initial data to be constant outside $\bar\Omega_0$ and focus on the effects of kernel functions and geometry property of initial domains.

\medskip

\noindent  \textbf{Example 1.} {\it The assumption \textbf{(K1)} on  kernel functions.}

For simplicity, we focus on the the one dimensional case and assume that the initial  domain is an internal. According  to Theorem \ref{thm-continuity}, if   the kernel function $k(x)$ is  decreasing in $|x|$, then $\Omega(t)$ expands continuously.

Now we construct an example where the  kernel function  is a perturbation of the following function
$$
k_*(x) =
\begin{cases}
 \displaystyle  {1\over 4\sigma}  &  1-\sigma \leq |x| \leq 1+\sigma, \\
  0 & \textrm{otherwise},
 \end{cases}
 $$
where $\displaystyle 0< \sigma< {1\over 4}$ is small, and verify that the  jumping phenomena   happens for proper choice of initial data.

To be more specific, consider the problem
\begin{equation}\label{Ex1-approx}
\begin{cases}
\displaystyle  \gamma_t(t,x)=\int_{\mathbb R} k_j(x-y)\gamma^+(t,y)dy-\gamma^+(t,x)   &  t>0,\ x\in\mathbb R, \\
\displaystyle\gamma(0,x)=c_0   & x\in \left(-{1\over 4},{1\over 4}\right) ,\\
\displaystyle\gamma(0,x)=-\ell_0   & x\in \mathbb R  \setminus \left(-{1\over 4},{1\over 4}\right),
\end{cases}
\end{equation}
where $c_0$, $\ell_0$ are positive constants,  $k_j$ satisfies the assumption \textbf{(K)} and
$$
\lim_{j\rightarrow \infty}\| k_j -k_*\|_{L^1(\mathbb R^n)}=0.
$$
Let $\gamma_j$ denote the solution to the problem (\ref{Ex1-approx}).

\medskip

We claim that {\it if  $\displaystyle 2 \ell_0 <   c_0$, $\displaystyle 0<\sigma<{1\over 4}$, then  the jumping phenomena happens for  the nonlocal Stefan problem (\ref{Ex1-approx}) when $j$ is sufficiently large.}

\medskip

To prove the cliam, first consider the problem
\begin{equation}\label{Ex1-0}
\begin{cases}
\displaystyle  \gamma_t(t,x)=\int_{\mathbb R} k_*(x-y)\gamma^+(t,y)dy- \gamma^+(t,x)   &  t>0,\ x\in\mathbb R, \\
\displaystyle\gamma(0,x)=c_0    & x\in \left(-{1\over 4},{1\over 4}\right) ,\\
\displaystyle\gamma(0,x)=-\ell_0   & x\in \mathbb R  \setminus \left(-{1\over 4},{1\over 4}\right).
\end{cases}
\end{equation}
The existence and uniqueness of the solution, denoted by $\gamma_*$, to this problem  can be verified by similar arguments to the proof of Theorem \ref{thm-wellposedness}.
Moreover, similar to the proof of (\ref{pf-thm-ex1}) in the proof of Theorem \ref{thm-wellposedness}, one has
$$
\lim_{j\rightarrow \infty}\| \gamma_j -\gamma_*\|_{L^{\infty}(\mathbb R^n)}=0.
$$
Hence it suffices to show that the jumping phenomena happens in the limiting problem (\ref{Ex1-0})  if  $\displaystyle 2 \ell_0 <   c_0$, $\displaystyle 0<\sigma<{1\over 4}$.

Let $t_1$  denote the moment when $\gamma_*$ first  touches zero somewhere in  $\displaystyle \mathbb R  \setminus \left(-{1\over 4},{1\over 4}\right)$.
For  $\displaystyle x\in \left(-{1\over 4}, {1\over 4}\right)$, $0<t<t_1$, it is easy to see that $\displaystyle\int_{\mathbb R^n} k_*(x-y)\gamma_*^+(t,y)dy =0$ due to the definition of $k_*$ and the choice of $\sigma$. Thus
$$
\begin{cases}
(\gamma_*^+)_t(t,x)= -\gamma_*^+(t,x)   &0<t<t_1,\  x\in \left(-{1\over 4}, {1\over 4}\right),\\
\gamma_*^+(0,x)=c_0   & x\in \left(-{1\over 4},{1\over 4}\right),
\end{cases}
$$
thus
$$
\gamma_*^+(t,x) = c_0 e^{-t}\ \ \textrm{for}\  0<t<t_1,\ x\in  \left(-{1\over 4}, {1\over 4}\right).
$$

Now, for any fixed  $\displaystyle x^* \in \left\{ x\in \mathbb R  \setminus \left(-{1\over 4},{1\over 4}\right) \ | \ \gamma_*(t_1, x)=0 \right\}  $, we compute
\begin{eqnarray}\label{pf-l0}
\ell_0 = \int_0^{t_1} \int_{-{1\over 4}}^{{1\over 4}}k_* (x^*-y)c_0 e^{-t} dy dt = c_0 \left(1- e^{-t_1} \right)\int_{-{1\over 4}}^{{1\over 4}}k_* (x^*-y)  dy.
\end{eqnarray}
According to the definition of $k_*$, it is routine to verify that
$\displaystyle\int_{-{1\over 4}}^{{1\over 4}}k_* (x^*-y)  dy\leq {1\over 2}$ and
$$
\int_{-{1\over 4}}^{{1\over 4}}k_* (x^*-y)  dy = {1\over 2}\  \ \textrm{if and only if}\ \  x^* \in \left[-{5\over 4}+\sigma, -{3\over 4}-\sigma\right] \bigcup \left[{3\over 4}+ \sigma, {5\over 4}- \sigma\right].
$$
Hence when $ 2  \ell_0<  c_0$, one has
$$
\left\{ x\in \mathbb R  \setminus \left(-{1\over 4},{1\over 4}\right) \ \Big| \ \gamma(t_1, x)=0 \right\} = \left[-{5\over 4}+\sigma, -{3\over 4}-\sigma\right] \bigcup \left[{3\over 4}+ \sigma, {5\over 4}- \sigma\right],
$$
where
$$
t_1= -\ln \left( 1- {2\ell_0\over c_0}\right),
$$
i.e.
$$
 \{ x\in\mathbb R^n \ |\ \gamma_*(t_1, x)\geq 0 \} =\left(-{1\over 4}, {1\over 4}\right)\bigcup \left[-{5\over 4}+\sigma, -{3\over 4}-\sigma\right] \bigcup \left[{3\over 4}+ \sigma, {5\over 4}- \sigma\right]
$$
and  for $0\leq t<t_1$,
$$
\{ x\in\mathbb R^n \ |\ \gamma_*(t, x)\geq 0 \} =\left(-{1\over 4}, {1\over 4}\right).
$$

Therefore, the desired jumping phenomena happens in the problem (\ref{Ex1-0}) provided that $\displaystyle 0<\sigma<{1\over 4}$ and $\displaystyle 2  \ell_0<   c_0$.

\medskip

\noindent  \textbf{Example 2.} {\it The condition on the shape of initial domain $\bar\Omega_0$.}

To emphasize the effect of geometry property of initial domains, we still assume that   the kernel functions constructed in this example   satisfy the conditions \textbf{(K)} and \textbf{(K1)} as required in Theorem \ref{thm-continuity}. Then according  to Theorem \ref{thm-continuity},  if the initial domain is convex, $\Omega(t)$ expands continuously. However, in the   following  example,  we demonstrate that  the    jumping phenomena could happen when the initial domain is non-convex.

Consider the problem
\begin{equation}\label{Ex2-approx}
\begin{cases}
\displaystyle \gamma_t(t,x)=\int_{\mathbb R^n} \tilde k_j(x-y)\gamma^+(t,y)dy-\gamma^+(t,x)   & t>0,\  x\in\mathbb{R}^n,\\
\gamma(0,x)=c_0    & x\in \bar\Omega_0 ,\\
\gamma(0,x)=-\ell_0    & x\in \mathbb{R}^n \setminus \bar\Omega_0,
\end{cases}
\end{equation}
where $c_0$, $\ell_0$ are positive constants, $n\geq 2$, $\bar\Omega_0 = \bar B_2(0)\setminus B_1(0)$,  the kernel function $\tilde k_j$ satisfies the conditions \textbf{(K)} and \textbf{(K1)} for kernel functions as required in Theorem \ref{thm-continuity} and
$$
\lim_{j\rightarrow \infty}\| \tilde k_j -\tilde k\|_{L^1(\mathbb R^n)}=0,
$$
where $\omega_n$ denotes the volume of a unit ball in $\mathbb R^n$ and
$$
\tilde k(x) =
\begin{cases}
 \displaystyle  2^{-n} \omega_n^{-1}  &   |x| \leq 2, \\
  0 & \textrm{otherwise}.
 \end{cases}
 $$
The existence of such kernel functions  $\tilde k_j$  is obvious.

\medskip
We claim that {\it if $\displaystyle \left(1-2^{-n} \right)c_0 > \ell_0$, then for  $j$ sufficiently  large,  the jumping phenomena happens for the nonlocal Stefan problem  (\ref{Ex2-approx}).}
\medskip

Similar to  \textbf{Example 1}, to prove this claim, it suffices to show the jumping phenomena happens for the following model
\begin{equation}\label{Ex2-0}
\begin{cases}
\displaystyle \gamma_t(t,x)=\int_{\mathbb R^n} \tilde k(x-y)\gamma^+(t,y)dy- \gamma^+(t,x)   & t>0,\ x\in\mathbb{R}^n, \\
\gamma(0,x)=c_0     & x\in \bar\Omega_0 ,\\
\gamma(0,x)=-\ell_0    & x\in \mathbb{R}^n \setminus \bar\Omega_0,
\end{cases}
\end{equation}
if $\displaystyle \left(1-2^{-n} \right)c_0 > \ell_0$.

Let $\tilde\gamma $ denote the solution to the problem (\ref{Ex2-0}) and $t_2$, if exists,  denote the moment when  the solution $\tilde\gamma$   first touches zero somewhere in  $\mathbb{R}^n \setminus \bar\Omega_0$.
When $0<t< t_2$, for any $x\in \mathbb{R}^n \setminus \bar\Omega_0$ with $x\neq 0$, thanks to the definition of $\tilde k$, it is easy to check that
\begin{eqnarray*}
\tilde\gamma (t,x) + \ell_0& =&  \int_0^t\int_{\mathbb R^n} \tilde k(x-y)\tilde\gamma^+(\tau,y)dy d\tau = \int_0^t \int_{\bar B_2(0)\setminus B_1(0)} \tilde k(x-y)\tilde\gamma^+(\tau,y)dy d\tau\\
&<&  \int_0^t \int_{\bar B_2(0)\setminus B_1(0)} \tilde k(-y)\tilde\gamma^+(\tau,y)dy d\tau=\tilde\gamma (t,0) + \ell_0.
\end{eqnarray*}
This indicates that if $t_2<+\infty$, then at $t=t_2$, $\tilde\gamma$ touches zero only at $x=0$, i.e.,
$$
 \{ x\in\mathbb R^n \ |\ \tilde\gamma(t_1, x)\geq 0 \} =\bar B_2(0)\setminus B_1(0) \ \bigcup\ \{x=0\}
$$
and  for $0\leq t<t_2$,
$$
\{ x\in\mathbb R^n \ |\ \tilde\gamma(t, x)\geq 0 \} =\bar B_2(0)\setminus B_1(0).
$$
The jumping phenomena happens.

It remains to show that $t_2<+\infty$. Suppose that $t_2 =+\infty$. Based on the definition of $\tilde k$ and the first equation in (\ref{Ex2-0}), it is easy to see that
\begin{eqnarray*}
\ell_0 \geq  \int_0^{+\infty}  \int_{\bar B_2(0)\setminus B_1(0)} \tilde k(-y) \tilde\gamma^+(t,y) dy dt > \int_0^{+\infty} \left(1-2^{-n} \right)c_0 e^{-t} dt =\left(1-2^{-n} \right)c_0 > \ell_0.
\end{eqnarray*}
This is impossible.
The proof is complete.

\section{Proofs of main results}

\subsection{Equivalent characterization of optimal convergence condition}
This section is devoted to the proof of Proposition \ref{prop-kernel}, where we derive an equivalent characterization of optimal convergence condition on the kernel functions.

\begin{proof}[Proof of Proposition \ref{prop-kernel}]
Assume that (i) holds. For clarity, set $\displaystyle w= (w_1, ..., w_n) ={\xi\over |\xi|}$.
Then we compute as follows.
\begin{eqnarray}
\frac{1-\hat k (\xi) }{|\xi|^2} &=& {1\over |\xi|^2} \left( 1- \int_{\mathbb R^n} e^{-ix\cdot \xi} k(x)dx \right)\cr
&=&  {1\over |\xi|^2}  \int_{\mathbb R^n}  ix\cdot w   \int_0^{|\xi|} e^{- (ix\cdot w)   \eta}d\eta k(x)dx\cr
&=&  {1\over |\xi|^2}  \int_{\mathbb R^n}   ix\cdot w  \int_0^{|\xi|} \left(  e^{-(ix\cdot w)   \eta} -1 \right)d\eta k(x)dx\cr
&=&  {1\over |\xi|^2}  \int_{\mathbb R^n}   (x\cdot w)^2  \int_0^{|\xi|} \int_0^{\eta}   e^{- (ix\cdot w)   \tau} d\tau  d\eta k(x)dx,\nonumber
\end{eqnarray}
where the third equality is due to the first equality in (i). Notice that thanks to the assumptiones in (i), we have
$$
 \int_{\mathbb R^n}   (x\cdot w)^2 k(x)dx ={1\over n}\int_{\mathbb R^n} |x|^2 k(x) dx.
$$
Then it follows that
\begin{eqnarray*}
&&\frac{1-\hat k (\xi) }{|\xi|^2} - {1\over 2n}\int_{\mathbb R^n} |x|^2 k(x) dx \\
&=&{1\over |\xi|^2}  \int_{\mathbb R^n}   (x\cdot w)^2  \int_0^{|\xi|} \int_0^{\eta}   e^{- (ix\cdot w)   \tau} d\tau  d\eta k(x)dx- {1\over 2} \int_{\mathbb R^n}   (x\cdot w)^2 k(x)dx\\
&= &{1\over |\xi|^2}  \int_{\mathbb R^n}   (x\cdot w)^2  \int_0^{|\xi|} \int_0^{\eta}   \left(e^{- (ix\cdot w)   \tau}-1\right) d\tau  d\eta k(x)dx.
\end{eqnarray*}
Lebesgue dominated convergence theorem yields that
$$
\lim_{\xi \rightarrow 0} \frac{1-\hat k (\xi) }{|\xi|^2}- {1\over 2n}\int_{\mathbb R^n} |x|^2 k(x) dx =0.
$$
Thus (ii) is verified and $\displaystyle {1\over 2n}\int_{\mathbb R^n} |x|^2 k(x) dx =A.$

\bigskip

Assume that (ii) is valid. First choose $\xi= (0,..., \xi_j, ...,0)$, $1\leq j\leq n$, with $\xi_j>0$, then
\begin{eqnarray}\label{pf-lem-real-imagine}
&&\frac{1-\hat k (\xi) }{|\xi|^2} = {1\over |\xi|^2} \left( 1- \int_{\mathbb R^n} e^{-ix\cdot \xi} k(x)dx \right)={1\over \xi_j^2}  \int_{\mathbb R^n} \left( 1- e^{-i x_j \xi_j} \right) k(x)dx\cr
&=& {1\over \xi_j^2}  \int_{\mathbb R^n} \left( 1- \cos (x_j \xi_j) + i\sin  (x_j \xi_j)\right) k(x)dx.
\end{eqnarray}
For any $R>0$, we have
$$
\frac{| 1-\hat k (\xi)| }{|\xi|^2} \geq  {1\over \xi_j^2}  \int_{ B_R(0)}  \left(  1- \cos (x_j \xi_j) \right)k(x)dx,
$$
which yields that
$$
\lim_{\xi_j \rightarrow 0} \frac{|1-\hat k (\xi)| }{|\xi|^2} \geq {1\over 2}  \int_{B_R(0)} x_j^2 k(x)dx.
$$
Since $R$ is arbitrary and $1\leq j\leq n$, one sees that
\begin{equation}\label{pf-lem-x-2}
   {1\over 2n}\int_{\mathbb R^n} |x|^2 k(x) dx\leq A<+\infty.
\end{equation}
This also indicates that
\begin{equation}\label{pf-lem-x-1}
\displaystyle \int_{\mathbb R^n} |x| k(x) dx<+\infty.
\end{equation}

Next, we still choose $\xi= (0,..., \xi_j, ...,0)$, $1\leq j\leq n$, with $\xi_j>0$. Notice that
\begin{eqnarray*}
\frac{1-\hat k (\xi) }{|\xi|} = 
{1\over \xi_j}  \int_{\mathbb R^n} \left( 1- e^{-i x_j \xi_j} \right) k(x)dx =  {1\over \xi_j}  \int_{\mathbb R^n} ix_j \int_0^{\xi_j} e^{-i x_j \eta} d\eta k(x)dx,
\end{eqnarray*}
where $x=(0,...,x_j,...,0)$.
Due to (\ref{pf-lem-x-1}), Lebesgue dominated convergence theorem can be applied and  one sees that
$$
0=\lim_{\xi \rightarrow 0} \frac{1-\hat k (\xi) }{|\xi|} = \int_{\mathbb R^n} ix_j  k(x)dx,
$$
i.e.,
\begin{equation}\label{pf-lem-xk(x)-0}
\displaystyle \int_{\mathbb R^n} x_jk(x)dx=0,  \ 1\leq j\leq n.
\end{equation}

Now thanks to (\ref{pf-lem-xk(x)-0}), we have
\begin{eqnarray*}
&& {1\over \xi_j^2}  \int_{\mathbb R^n}\sin  (x_j \xi_j) k(x)dx\\
&=& {1\over \xi_j^2}  \int_{\mathbb R^n} x_j \int_0^{\xi_j}  \cos  (x_j \eta)  d\eta k(x)dx={1\over \xi_j^2}  \int_{\mathbb R^n} x_j \int_0^{\xi_j} \left( \cos  (x_j \eta) -1 \right) d\eta k(x)dx\\
&=& {1\over \xi_j^2}  \int_{\mathbb R^n} -x_j^2 \int_0^{\xi_j} \int_0^{\eta}\sin  (x_j \tau) d\tau d\eta k(x)dx.
\end{eqnarray*}
Thus (\ref{pf-lem-x-2}) and Lebesgue dominated convergence theorem  imply that
$$
\lim_{\xi_j \rightarrow 0}{1\over \xi_j^2}  \int_{\mathbb R^n}\sin  (x_j \xi_j) k(x)dx =0.
$$
Now in (\ref{pf-lem-real-imagine}), letting $\xi_j \rightarrow 0$, again it follows from  (\ref{pf-lem-x-2}), (\ref{pf-lem-xk(x)-0}) and  Lebesgue dominated convergence theorem that
\begin{eqnarray*}
A &= &\lim_{\xi_j \rightarrow 0} {1\over \xi_j^2}  \int_{\mathbb R^n} \left( 1- \cos (x_j \xi_j) \right) k(x)dx  = \lim_{\xi_j \rightarrow 0} {1\over \xi_j^2}  \int_{\mathbb R^n} x_j\int_0^{\xi_j} \sin (x_j \eta) d\eta k(x)dx\\
&=& \lim_{\xi_j \rightarrow 0} {1\over \xi_j^2}  \int_{\mathbb R^n} x_j\int_0^{\xi_j} \left(\sin (x_j \eta)-1\right) d\eta k(x)dx\\
&=& \lim_{\xi_j \rightarrow 0} {1\over \xi_j^2}  \int_{\mathbb R^n} x_j^2\int_0^{\xi_j} \int_0^{\eta} \cos (x_j\tau)d\tau d\eta k(x)dx = {1\over 2} \int_{\mathbb R^n} x_j^2  k(x)dx.
\end{eqnarray*}
Hence
\begin{equation}\label{pf-lem-x2k(x)-equal}
\int_{\mathbb R^n} x_j^2 k(x)dx =2A ={1\over n}\int_{\mathbb R^n} |x|^2 k(x) dx, \ \ 1\leq j\leq n.
\end{equation}

At the end, it remains to show that  $\displaystyle\int_{\mathbb R^n} x_jx_h k(x)dx =0$, $1\leq j,h\leq n, j\neq h.$ Choose $\xi= (0,..., \xi_j,..., \xi_h, ...,0)$ with $j<h$, $\xi_j>0$ and $\xi_h= \lambda \xi_j$. Then thanks to (\ref{pf-lem-xk(x)-0})  and (\ref{pf-lem-x2k(x)-equal}), it follows that
\begin{eqnarray*}
&&\frac{1-\hat k (\xi) }{|\xi|^2} = {1\over |\xi|^2} \left( 1- \int_{\mathbb R^n} e^{-ix\cdot \xi} k(x)dx \right)={1\over \xi_j^2+ \lambda^2\xi_j^2}  \int_{\mathbb R^n} \left( 1- e^{-i \left(x_j +\lambda x_h \right)\xi_j} \right) k(x)dx\cr
&=&  {1\over \xi_j^2+ \lambda^2\xi_j^2}  \int_{\mathbb R^n}   i \left(x_j +\lambda x_h \right)
\int_0^{\xi_j} e^{-i \left(x_j +\lambda x_h \right)\eta} d \eta   k(x)dx\cr
&=&  {1\over \xi_j^2+ \lambda^2\xi_j^2}  \int_{\mathbb R^n}   i \left(x_j +\lambda x_h \right)
\int_0^{\xi_j} \left( e^{-i \left(x_j +\lambda x_h \right)\eta}-1\right) d \eta   k(x)dx\cr
&=&  {1\over \xi_j^2+ \lambda^2\xi_j^2}  \int_{\mathbb R^n}   \left(x_j +\lambda x_h \right)^2
\int_0^{\xi_j} \int_0^{\eta} e^{-i \left(x_j +\lambda x_h \right)\tau} d\tau d \eta   k(x)dx.
\end{eqnarray*}
Letting $\xi_j \rightarrow 0$, Lebesgue dominated convergence theorem and (\ref{pf-lem-x2k(x)-equal}) imply that
\begin{eqnarray*}
A= {1\over 2(1+ \lambda^2)}  \int_{\mathbb R^n}   \left(x_j +\lambda x_h \right)^2
   k(x)dx =A+ {\lambda \over 1+ \lambda^2}  \int_{\mathbb R^n}   x_j  x_h k(x)dx.
\end{eqnarray*}
This indicates that
$$
\int_{\mathbb R^n}   x_j  x_h k(x)dx =0,
$$
since $\lambda$ is an arbitrary constant.
The proof is complete.
\end{proof}

\subsection{The nonlocal to local convergence of   Stefan problem}
This section is devoted to the proof of Theorem \ref{thm-convergence}, where the convergence relations between local and nonlocal one-phase Stefan problems are verified under the optimal condition  (\ref{assumption-hat k}) imposed on the kernel function as follows
\begin{equation*}
\hat k (\xi )  = 1- A|\xi |^2 +o(|\xi |^2) \ \ \ \textrm{as}  \  \xi\rightarrow 0,
\end{equation*}
where $A>0$ is a constant.

For convenience,  recall the classical one-phase problem (\ref{localstephan-n}) with $d=A$ as follows
\begin{equation*}
\begin{cases}
\displaystyle \theta_t(t,x)= A \Delta\theta (t,x)  & t>0,\ x\in \{\theta(t, \cdot)>0\}, \\
\nabla_x \theta \cdot \nabla_x s= -\ell_0   &t>0,\ x\in\partial \{\theta(t, \cdot)>0\},  \\
\theta=0 &t>0,\ x\in\partial \{\theta(t, \cdot)>0\},  \\
\theta(0,x)=u_0(x)   & x\in \bar\Omega_0,
\end{cases}
\end{equation*}
where the free boundary $\partial \{\theta(t, \cdot)>0\}$ at the time $t$ is given by the equation $s(x)=t$ and set $s(x)=0$ if $x\in \bar\Omega_0$.
It is known that the classical one-phase problem can be reduced to a parabolic variational inequality \cite[Chapter 1.9]{Friedman}. To be more specific, define
\begin{equation*}
v(t,x) =
\begin{cases}
\displaystyle \int_{0}^t \theta(\tau, x)d\tau& \textrm{if}\ x\in \bar\Omega_0,\\
\displaystyle 0  & \textrm{if}\ x\in\mathbb R^n\setminus\bar\Omega_0,\, t\leq s(x),\\
\displaystyle \int_{s(x)}^t \theta(\tau, x)d\tau & \textrm{if}\ x\in\mathbb R^n\setminus\bar\Omega_0,\, t>s(x),
\end{cases}
\end{equation*}
and then  the problem (\ref{localstephan-n}) is transformed  into a   variational inequality for the function $v(t,x)$ as follows
\begin{equation}\label{variationalinequality}
\begin{cases}
v_t-A \Delta v\geq \bar f & \textrm{a.e. in}\ (0,T)\times \mathbb R^n,\\
v\geq 0   &  \textrm{a.e. in}\ (0,T)\times \mathbb R^n,\\
(v_t- A \Delta v- \bar f)v=0    &  \textrm{a.e. in}\ (0,T)\times \mathbb R^n,
\end{cases}
\end{equation}
where $\bar{f}=\gamma_0$ is defined in (\ref{initial-onephase}).  It has been proved that there exists a unique solution of the problem (\ref{variationalinequality}), still denoted by $v(t,x)$, and
$$
D_x v,\, D_x^2 v,\, D_t v \ \  \ \textrm{belong to}\ L^{\infty}((0,T); L^p(\mathbb R^n))\  \ \textrm{for}\ p<\infty.
$$
See \cite[Chapter 1.9]{Friedman} for details.

Borrowing this idea, for the nonlocal one-phase Stefan problem (\ref{one-nonlocalstefan-epsilon}),
define
\begin{equation}\label{v-epsilon}
\displaystyle v_{\epsilon}(t,x) = \int_0^t \gamma^+_{\epsilon}(\tau,x)  d\tau,
\end{equation}
Theorem \ref{thm-convergence} is about the convergence relations between $v_{\epsilon}$ defined in (\ref{v-epsilon}) and the solution $v$ to the problem  (\ref{variationalinequality}).

First of all, we compute the equation satisfied by  $v_{\epsilon}$. For any $x\in \mathbb R ^n \setminus \bar\Omega_0$, let $s_{\epsilon}(x)$ denote the time if exists when $\gamma_{\epsilon}(t,x)$  first reaches zero. Thus
\begin{equation}\label{boundary}
\ell_0   ={1\over \epsilon^2} \int_0^{s_{\epsilon}(x)} \int_{\mathbb R ^n}k_{\epsilon}(x-y) \gamma_{\epsilon}^+(\tau ,y) dy d\tau.
\end{equation}
\begin{itemize}
\item if $x\in \bar\Omega_0,  t>0$, then
\begin{eqnarray}
&& v_{\epsilon t} - {1\over \epsilon^2}  \int_{\mathbb R^n} k_{\epsilon}(x-y)v_{\epsilon}(t,y)dy  + {1\over \epsilon^2}  v_{\epsilon}(t,x) \cr
&=& \gamma_{\epsilon}^+(t,x) - {1\over \epsilon^2}  \int_{\mathbb R^n} k_{\epsilon}(x-y)\int_0^t \gamma_{\epsilon}^+(\tau, y) d\tau dy  + {1\over \epsilon^2}  \int_0^t \gamma_{\epsilon}^+(\tau, x) d\tau\cr
&=& \int_0^t \gamma_{\epsilon t}^+(\tau, x)d\tau + \gamma_{\epsilon}^+(0,x) - {1\over \epsilon^2} \int_{\mathbb R^n} k_{\epsilon}(x-y)\int_0^t \gamma_{\epsilon}^+(\tau, y) d\tau dy  +{1\over \epsilon^2}  \int_0^t \gamma_{\epsilon}^+(\tau, x) d\tau\cr
&=&  \gamma_0(x);\nonumber
\end{eqnarray}
\item  if $x\in \mathbb R ^n \setminus \bar\Omega_0,  0<t\leq s_{\epsilon}(x)$, then $v_{\epsilon}(t,x)=0$. Thus
\begin{eqnarray*}
&& v_{\epsilon t} - {1\over \epsilon^2}  \int_{\mathbb R^n} k_{\epsilon}(x-y)v_{\epsilon}(t,y)dy  + {1\over \epsilon^2} v_{\epsilon}(t,x) = -  {1\over \epsilon^2}  \int_{\mathbb R^n} k_{\epsilon}(x-y)v_{\epsilon}(t,y)dy;
\end{eqnarray*}
\item  if $x\in \mathbb R ^n \setminus \bar\Omega_0,  t > s_{\epsilon}(x)$, then
\begin{eqnarray}
&& v_{\epsilon t} - {1\over \epsilon^2}  \int_{\mathbb R^n} k_{\epsilon}(x-y)v_{\epsilon}(t,y)dy  + {1\over \epsilon^2} v_{\epsilon}(t,x) \cr
&=& \gamma_{\epsilon}^+(t,x) - {1\over \epsilon^2} \int_{\mathbb R^n} k_{\epsilon}(x-y)\int_0^t \gamma_{\epsilon}^+(\tau, y) d\tau dy  +{1\over \epsilon^2}  \int_0^t \gamma_{\epsilon}^+(\tau, x) d\tau\cr
&=& \int_{s_{\epsilon}(x)}^t \gamma_{\epsilon t}^+(\tau, x)d\tau  - {1\over \epsilon^2} \int_{\mathbb R^n} k_{\epsilon}(x-y)\int_{s_{\epsilon}(x)}^t \gamma_{\epsilon}^+(\tau, y) d\tau dy  + {1\over \epsilon^2} \int_{s_{\epsilon}(x)}^t \gamma_{\epsilon}^+(\tau, x) d\tau \cr
&&-{1\over \epsilon^2}  \int_{\mathbb R^n} k_{\epsilon}(x-y)\int_0^{s_{\epsilon}(x)} \gamma_{\epsilon}^+(\tau, y) d\tau dy\cr
&=& -\ell_0,\nonumber
\end{eqnarray}
according to (\ref{boundary}).
\end{itemize}
Hence one sees that $v_{\epsilon}$ satisfies
\begin{equation}\label{nonlocalstephan-n-v-epsilon}
\begin{cases}
\displaystyle v_{\epsilon t}(t,x)= {1\over \epsilon^2} \int_{\mathbb R^n} k_{\epsilon}(x-y)v_{\epsilon}(t,y)dy- {1\over \epsilon^2} v_{\epsilon}(t,x) + f_{\epsilon}(t,x )    & t>0,\ x\in \mathbb R^n,\\
v_{\epsilon}(0,x)=0    & x\in \mathbb{R}^n,
\end{cases}
\end{equation}
where
$$
 f_{\epsilon}(t,x)=
\begin{cases}
\displaystyle \gamma_0(x)   & t>0,\ x\in \bar\Omega_0, \\
\displaystyle -  \int_{\mathbb R^n}{1\over \epsilon^2} k_{\epsilon}(x-y)v_{\epsilon}(t,y)dy  & 0< t\leq s_{\epsilon}(x),\ x\in \mathbb R^n \setminus \bar\Omega_0,\\
\displaystyle -\ell_0  & t>s_{\epsilon}(x),\  x\in \mathbb R^n \setminus \bar\Omega_0.
\end{cases}
$$

Secondly, we prepare  some useful estimates about    $f_{\epsilon}$.

\begin{lem}\label{lem-f-epsilon-estimate}
Assume  that   in the problem (\ref{one-nonlocalstefan-epsilon}), the kernel  function  $k$ satisfies the assumption  \textbf{(K)}, the condition (\ref{initialdata}) is valid and  the initial data satisfies (\ref{initial-onephase}).  Then $f_\epsilon(t,x)$ is uniformly bounded in $L^p(\mathbb{R}^n)$ for any $\epsilon>0$, $t>0$,  $1\leq p\leq \infty$.
\end{lem}



\begin{proof}
Similar to the proof of Proposition \ref{prop-Lp-estimate},  if $\phi\in L^1(\mathbb R^n)\bigcap L^{\infty}(\mathbb R^n)$, then  for any $p>1$, $\phi \in L^p(\mathbb R^n)$ and
$$
\|\phi\|_{L^p(\mathbb R^n)}\leq \left(\|\phi\|^{p-1}_{L^{\infty}(\mathbb R^n)}\|\phi\|_{L^1(\mathbb R^n)} \right)^{1\over p}\leq \left(\|\phi\|_{L^{\infty}(\mathbb R^n)}+1\right)\left(\|\phi\|_{L^1(\mathbb R^n)}+1 \right).
$$
Hence it suffices to verify the conclusion  for $p=1$ and $p=\infty$.

 Since  $f_\epsilon(t,x)= \gamma_0$ for $x\in\bar{\Omega}_0$ and $t>0$,  we only need to estimate $f_\epsilon$ outside $\bar{\Omega}_0$.
It mainly relies on the following estimates:
\begin{equation}\label{pf-lem-key}
-f_\epsilon(t,x) \in [0, \ell_0]    \ \ \textrm{for}\  x\in \mathbb R^n \setminus \bar\Omega_0,
 \ \ \int_{\mathbb{R}^n\setminus \bar{\Omega}_0}-f_\epsilon(t,x)dx\leq \int_{\Omega_0}\gamma_0dx.
\end{equation}

Assume that  (\ref{pf-lem-key})  is valid. It immediately yields that
$$
\| f_{\epsilon}(t, \cdot)\|_{L^{\infty}(\mathbb R^n)}\leq \max\, \left\{\|\gamma_0|_{\bar\Omega_0}\|_{L^{\infty}(\Omega_0)}, \ell_0 \right\},\ \
\| f_{\epsilon}(t, \cdot)\|_{L^1(\mathbb R^n)}\leq 2\int_{\Omega_0}\gamma_0dx.
$$
The desired conclusion follows.

Now it remains to verify (\ref{pf-lem-key}). In fact,  due to (\ref{boundary}), the first estimate in (\ref{pf-lem-key}) is obvious.
Intuitively, the second estimate in (\ref{pf-lem-key}) indicates that    $\displaystyle\int_{\mathbb{R}^n \setminus\bar{\Omega}_0} -f_\epsilon(t,x)dx$ is less than  the total energy absorbed outside $\bar\Omega_0$ from time 0 to $t$, which can not exceed the total energy at the initial time, i.e. $\displaystyle\int_{\Omega_0}\gamma_0dx$. To be more precise, by (\ref{one-nonlocalstefan-epsilon}), one has for any large $R>0$
\begin{eqnarray}\label{lem-ii}
&&\int_{B_R(0)\setminus \bar\Omega_0}\left(\gamma_{\epsilon }(t,x) -\gamma_{\epsilon }(0,x)  \right)dx\cr
&=& \int_{(B_R(0)\setminus \bar\Omega_0 ) \bigcap \{ s_{\epsilon}(x)<t \}}\left(\gamma_{\epsilon }(t,x) -\gamma_{\epsilon }(0,x)  \right)dx +\int_{(B_R(0)\setminus \bar\Omega_0)  \bigcap \{ s_{\epsilon}(x)\geq t \}}\left(\gamma_{\epsilon }(t,x) -\gamma_{\epsilon }(0,x)  \right)dx\cr
&\geq &\int_{(B_R(0)\setminus \bar\Omega_0 ) \bigcap \{ s_{\epsilon}(x)<t \}} \ell_0dx + {1\over \epsilon^2}\int_0^t \int_{(B_R(0)\setminus \bar\Omega_0)  \bigcap \{ s_{\epsilon}(x)\geq t \}} \int_{\mathbb R^n} k_{\epsilon}(x-y)\gamma_{\epsilon}^+(\tau,y)dy dx d\tau\cr
&=& \int_{(B_R(0)\setminus \bar\Omega_0 ) \bigcap \{ s_{\epsilon}(x)<t \}} \ell_0dx + {1\over \epsilon^2} \int_{(B_R(0)\setminus \bar\Omega_0)  \bigcap \{ s_{\epsilon}(x)\geq t \}} \int_{\mathbb R^n} k_{\epsilon}(x-y) v_{\epsilon}(t,y)dy dx \cr
&=& \int_{B_R(0)\setminus \bar\Omega_0} -f_\epsilon(t,x)dx.
\end{eqnarray}
Moreover, it is easy to see that
\begin{eqnarray}
&& \int_{B_R(0)  }\gamma_{\epsilon t}(t,x) dx = {1\over \epsilon^2}\int_{B_R(0)  } \int_{\mathbb R^n} k_{\epsilon}(x-y)\gamma_{\epsilon}^+(t,y)dy dx- {1\over \epsilon^2}\int_{B_R(0)  } \gamma_{\epsilon}^+(t,x)dx\cr
&\leq &  {1\over \epsilon^2} \int_{\mathbb R^n}  \gamma_{\epsilon}^+(t,x)dx- {1\over \epsilon^2}\int_{B_R(0)  } \gamma_{\epsilon}^+(t,x)dx,\nonumber
\end{eqnarray}
where the validity of the above inequality is due to the property that  $\gamma_{\epsilon}^+(t,\cdot)\in L^1(\mathbb R^n)$ proved Propsition \ref{prop-Lp-estimate}.
This implies that
\begin{eqnarray}
&& \int_{B_R(0)\setminus \bar\Omega_0}\left(\gamma_{\epsilon }(t,x) -\gamma_{\epsilon }(0,x)  \right)dx \cr
&\leq&  -\int_{ \bar\Omega_0}\left(\gamma_{\epsilon }(t,x) -\gamma_{\epsilon }(0,x)  \right)dx +{1\over \epsilon^2} \int_0^t \left(   \int_{\mathbb R^n}  \gamma_{\epsilon}^+(\tau,x)dx- \int_{B_R(0)  } \gamma_{\epsilon}^+(\tau,x)dx\right) d\tau\cr
&\leq&  \int_{\Omega_0}\gamma_0dx + {1\over \epsilon^2} \int_0^t \left(   \int_{\mathbb R^n}  \gamma_{\epsilon}^+(\tau,x)dx- \int_{B_R(0)  } \gamma_{\epsilon}^+(\tau,x)dx\right) d\tau.\nonumber
\end{eqnarray}
This, together with (\ref{lem-ii}), yields that
$$
\int_{B_R(0)\setminus \bar\Omega_0} -f_\epsilon(t,x)dx\leq \int_{\Omega_0}\gamma_0dx + {1\over \epsilon^2} \int_0^t \left(   \int_{\mathbb R^n}  \gamma_{\epsilon}^+(\tau,x)dx- \int_{B_R(0)  } \gamma_{\epsilon}^+(\tau,x)dx\right) d\tau,
$$
and thus the  second estimate  in (\ref{pf-lem-key}) follows immediately by letting $R\rightarrow \infty$.
\end{proof}


Moreover,  on the basis of the apriori estimates derived in Propositions \ref{prop-Lp-estimate} and \ref{prop-uniform-difference}, we establish some convergence results about $v_{\epsilon}$ defined in (\ref{v-epsilon}).

\begin{lem}\label{lem-convergence-ae}
Assume  that   in the problem (\ref{one-nonlocalstefan-epsilon}), the assumptions in Lemma \ref{lem-f-epsilon-estimate} are valid.  Then
for any fixed $t>0$, there exist a sequence $\{ \epsilon_\ell\}$, which depends on $t$ and satisfies $\lim_{\ell\rightarrow \infty} \epsilon_\ell =0$, and $\tilde{v}^t\in L^1(\mathbb R^n)$ such that
$v_{\epsilon_\ell}(t,\cdot) \rightarrow \tilde{v}^t(\cdot)$ a.e. in $\mathbb R^n$.
\end{lem}

\begin{proof}
Thanks to Proposition \ref{prop-uniform-difference},
\begin{eqnarray}
&& \int_{\mathbb R^n}\left|v_{\epsilon }(t,x+h) - v _{\epsilon}(t,x)\right|dx = \int_{\mathbb R^n}\left|\int_0^t \left( \gamma_{\epsilon}^+(\tau,x+h) -  \gamma_{\epsilon}^+(\tau,x) \right) d\tau \right|dx \cr
&\leq& \int_0^t \int_{\mathbb R^n} \left| \gamma_{\epsilon }(\tau,x+h) -  \gamma_{\epsilon}(\tau,x) \right|  dx d\tau  \cr
&\leq&  \int_0^t \int_{\mathbb R^n} \left|\gamma_0(x+h) -  \gamma_0(x) \right|  dx d\tau =t \int_{\mathbb R^n} \left|\gamma_0(x+h) -  \gamma_0(x) \right|  dx.\nonumber
\end{eqnarray}
This, together with the Fr$\acute{e}$chet-Kolmogorov theorem and Proposition \ref{prop-Lp-estimate}, indicates that for any fixed $t>0$ and  bounded set $\Omega\subseteq \mathbb R^n$, $\{v_{\epsilon}(t,\cdot) \, |\,  0<\epsilon<1 \}$ is precompact in $L^1(\Omega)$. Then  it is easy to show that there exist a sequence $\{ \epsilon_\ell\}$ with $\lim_{\ell\rightarrow \infty} \epsilon_\ell =0$ and $\tilde{v}^t\in L^1(\mathbb R^n)$ such that $v_{\epsilon_\ell}(t,\cdot) \rightarrow \tilde{v}^t(\cdot)$ in $L^1_{loc}(\mathbb R^n)$ and thus
$v_{\epsilon_\ell}(t,\cdot) \rightarrow \tilde{v}^t(\cdot)$ a.e. in $\mathbb R^n$.
\end{proof}

We emphasize that so far the optimal convergence condition (\ref{assumption-hat k}) has not been required yet.
After previous preparations, we are ready to complete the proof of Theorem \ref{thm-convergence}.

\begin{proof}[Proof of Theorem \ref{thm-convergence}]
From now on, fix $T>0$.
Back to the problem (\ref{nonlocalstephan-n-v-epsilon}) satisfied by $v_{\epsilon}$, by  the  Fourier transform and the property $\hat k_{\epsilon}(\xi) = \hat k (\epsilon\xi)$, we derive that
\begin{equation}\label{FT-v-epsilon}
\hat v_{\epsilon}(t, \xi) = \int_0^t  e^  {  {1\over \epsilon^2} \left( \hat k (\epsilon\xi) -1 \right)(t-\tau)  } \hat f_{\epsilon}(\tau, \xi)  d\tau.
\end{equation}
Due to the Parseval formula,
$$
\| f_{\epsilon}(t, \cdot)\|_{L^2(\mathbb R^n)} = \| \hat f_{\epsilon}(t, \cdot)\|_{L^2(\mathbb R^n)}.
$$
Then thanks to Lemma \ref{lem-f-epsilon-estimate},  there exist a sequence $\{ \epsilon_j\}$ with $\lim_{j\rightarrow \infty} \epsilon_j =0$ and $f_0,\, G_0   \in L^2((0,T) \times \mathbb R^n) $ such that
\begin{equation}\label{pf-f-weakconvergence}
\lim_{j\rightarrow \infty} f_{\epsilon_j}   = f_0 \ \ \textrm{weakly in }\ L^2((0,T) \times \mathbb R^n).
\end{equation}
and
\begin{equation}\label{pf-f-Fourier-weakconvergence}
\lim_{j\rightarrow \infty} \hat f_{\epsilon_j}   = G_0 \ \ \textrm{weakly in }\ L^2((0,T) \times \mathbb R^n).
\end{equation}

Notice that for any  test function $\psi(t, \xi)\in C_c((0,T) \times \mathbb R^n)$, on the one side, due to (\ref{pf-f-weakconvergence}),
\begin{eqnarray}
&&\lim_{j\rightarrow \infty}\int_0^T\int_{\mathbb R^n}\hat f_{\epsilon_j}(t, \xi)\psi(t, \xi) d\xi dt\cr
 &=&\lim_{j\rightarrow \infty} \int_0^T\int_{\mathbb R^n}\left( \int_{\mathbb R^n} e^{-ix\cdot \xi}f_{\epsilon_j}(t, x)dx \right)\psi(t, \xi) d\xi dt\cr
 &=& \lim_{j\rightarrow \infty}\int_0^T\int_{\mathbb R^n}\left( \int_{\mathbb R^n} e^{-ix\cdot \xi}\psi(t, \xi) d\xi \right) f_{\epsilon_j}(t, x)dx dt\cr
 &=& \int_0^T\int_{\mathbb R^n}\left( \int_{\mathbb R^n} e^{-ix\cdot \xi}\psi(t, \xi) d\xi \right) f_0(t, x)dx dt = \int_0^T\int_{\mathbb R^n}\hat f_0(t, \xi)\psi(t, \xi) d\xi dt.\nonumber
\end{eqnarray}
On the other side,  (\ref{pf-f-Fourier-weakconvergence}) yields that
$$
\lim_{j\rightarrow \infty} \int_0^T\int_{\mathbb R^n}\hat f_{\epsilon_j}(t, \xi)\psi(t, \xi) d\xi dt=\int_0^T\int_{\mathbb R^n} G(t, \xi)\psi(t, \xi) d\xi dt.
$$
Hence
\begin{equation*}
G_0(t,\xi) = \hat f_0\ \  \textrm{a.e.  in}\  (0,T) \times \mathbb R^n,
\end{equation*}
i.e.
\begin{equation}\label{pf-thm-G0}
\lim_{j\rightarrow \infty} f_{\epsilon_j}   = f_0,\ \ \lim_{j\rightarrow \infty} \hat f_{\epsilon_j}   = \hat f_0 \ \ \textrm{weakly in }\ L^2((0,T) \times \mathbb R^n).
\end{equation}

Introduce the following problem
\begin{equation}\label{problem-f0}
\begin{cases}
v_t = A \Delta v + f_0  & 0<t\leq T,\ x\in\mathbb R^n, \\
v(0, x) =0    & x\in\mathbb R^n,
\end{cases}
\end{equation}
where $v_*$ denote the unique generalized  solution   in $V_2^{1,1/2}(\mathbb R^n\times [0,T])$ \cite[Chapter III.5]{1968}. By  applying the Fourier transform to the problem (\ref{problem-f0}), we derive that
\begin{equation*}\label{FT-v}
\hat v_* (t, \xi)  =  \int_0^t  e^{ - A|\xi|^2 (t-\tau) } \hat f_0 (\tau, \xi)d\tau.
\end{equation*}

Fix $t\in(0,T)$. For any given $\phi(\xi)\in C_c^{\infty}(\mathbb R^n)$,
\begin{eqnarray}%
&& \lim_{j\rightarrow \infty}\int_{\mathbb R^n}\hat v_{\epsilon_j}(t, \xi)\phi(\xi)d\xi = \lim_{j\rightarrow \infty}\int_{\mathbb R^n}\left( \int_0^t  e^  {  {1\over \epsilon_j^2} \left( \hat k (\epsilon_j\xi) -1 \right)(t-\tau)  } \hat f_{\epsilon_j}(\tau, \xi)  d\tau\right) \phi(\xi)d\xi\cr
&=& \lim_{j\rightarrow \infty}\int_{\mathbb R^n}  \int_0^t  \left( e^  {  {1\over \epsilon_j^2} \left( \hat k (\epsilon_j\xi) -1 \right)(t-\tau)  }- e^{-A|\xi|^2(t-\tau)}\right) \hat f_{\epsilon_j}(\tau, \xi)  \phi(\xi) d\tau d\xi\cr
&& + \lim_{j\rightarrow \infty}\int_{\mathbb R^n}  \int_0^t    e^{-A|\xi|^2(t-\tau)}  \hat f_{\epsilon_j}(\tau, \xi)   \phi(\xi)d\tau d\xi.\nonumber
\end{eqnarray}
Since $\|\hat f_{\epsilon_j}(\tau, \cdot)\|_{L^{\infty}(\mathbb R^n)}\leq \| f_{\epsilon_j}(\tau, \cdot)\|_{L^1(\mathbb R^n)}$,  due to Lemma \ref{lem-f-epsilon-estimate}, the assumption (\ref{assumption-hat k}) and (\ref{pf-thm-G0}), we have
\begin{equation}\label{pf-thm-hat-v-1}
\lim_{j\rightarrow \infty}\int_{\mathbb R^n}\hat v_{\epsilon_j}(t, \xi)\phi(\xi)d\xi= \int_{\mathbb R^n}  \int_0^t    e^{-A|\xi|^2(t-\tau)}  \hat f_0(\tau, \xi)   \phi(\xi)d\tau d\xi=  \int_{\mathbb R^n} \hat v_* (t, \xi) \phi(\xi) d\xi.
\end{equation}

Moreover, thanks to Proposition \ref{prop-Lp-estimate}, there exists a subsequence of  $\{ \epsilon_j\}$,  denoted by $\{ \epsilon_{j_{\ell}}\}$, and $v_0^t$ in $L^2 (\mathbb R^n)$, such that $v_{\epsilon_{j_{\ell}}}(t, \cdot)\rightharpoonup v_0^t(\cdot)$ in $L^2 (\mathbb R^n)$. Then
for any given $\phi(\xi)\in C_c^{\infty}(\mathbb R^n)$,
\begin{eqnarray}\label{pf-thm-hat-v-2}
&& \lim_{{j_{\ell}}\rightarrow \infty}\int_{\mathbb R^n}\hat v_{\epsilon_{j_{\ell}}}(t, \xi)\phi(\xi)d\xi\cr
&=&\lim_{{j_{\ell}}\rightarrow \infty}\int_{\mathbb R^n} \left( \int_{\mathbb R^n} e^{-ix\cdot \xi} v_{\epsilon_{j_{\ell}}}(t, x) dx \right) \phi(\xi)d\xi \cr
&=&\lim_{{j_{\ell}}\rightarrow \infty}\int_{\mathbb R^n} \left( \int_{\mathbb R^n} e^{-ix\cdot \xi}  \phi(\xi)d\xi  \right) v_{\epsilon_{j_{\ell}}}(t, x) dx =\int_{\mathbb R^n} \left( \int_{\mathbb R^n} e^{-ix\cdot \xi}  \phi(\xi)d\xi  \right) v_0^t(x) dx\cr
&=& \int_{\mathbb R^n} \left( \int_{\mathbb R^n} e^{-ix\cdot \xi}   v_0^t(x) dx \right) \phi(\xi)d\xi.
\end{eqnarray}
Now (\ref{pf-thm-hat-v-1}) and (\ref{pf-thm-hat-v-2}) implies that
$$
\hat v_* (t, \xi) = \int_{\mathbb R^n} e^{-ix\cdot \xi}   v_0^t(x) dx\  \ \textrm{a.e. in}\ \mathbb R^n.
$$
Thus $v_* (t, x) =  v_0^t(x) $ a.e. in $\mathbb R^n$. Since $v_* (t, x)$ is the unique solution to the problem (\ref{problem-f0}), it follows immediately that for any $0<t<T$,  $v_{\epsilon}(t, \cdot)\rightharpoonup v_*(t,\cdot)$ in $L^2 (\mathbb R^n)$ as $\epsilon\rightarrow 0$.
This, together with Lemma \ref{lem-convergence-ae}, implies that \begin{equation}\label{pf-thm-ae}
v_{\epsilon}(t,x) \rightarrow v_*(t,x) \ \ \textrm{a.e. in} \ (0,T)\times \mathbb R^n\ \  \textrm{as}\  \epsilon\rightarrow 0.
\end{equation}

To complete the proof of Theorem \ref{thm-convergence}, it remains to verify that $v_*$ satisfies the parabolic variational inequality (\ref{variationalinequality}) as follows.
\begin{equation*}
\begin{cases}
v_t-A\Delta v\geq \bar f & \textrm{a.e. in}\ (0,T)\times \mathbb R^n,\\
v\geq 0   &  \textrm{a.e. in}\ (0,T)\times \mathbb R^n,\\
(v_t-A\Delta v- \bar f)v=0    &  \textrm{a.e. in}\ (0,T)\times \mathbb R^n,
\end{cases}
\end{equation*}
where $\bar{f}=\gamma_0$ for $x\in \mathbb{R}^n$.
Obviously $v_*\geq 0$ satisfies the first two inequalities in (\ref{variationalinequality}) since  $v_\epsilon $ is always non-negative and $f_\epsilon\geq \bar{f}$ for all $t>0$ and $x\in \mathbb{R}^n$.
Moreover, thanks to Lemma \ref{lem-f-epsilon-estimate},  (\ref{pf-f-weakconvergence}) and the uniqueness of weak convergence, it is standard to show that  $f_0 \in L^p(\mathbb R^n\times [0,T])$ for any $p>1$. Then by parabolic regularity theory and Sobolev embedding theorem, one obtains that $v_*(t,\cdot)$ is continuous in $\mathbb R^n$.
Thus, the set $\{v_* >0\}$ is open in $(0,T)\times \mathbb R^n$. Also notice that $f_\epsilon = \bar{f}$ if $v_{\epsilon}>0$. Hence thanks to (\ref{pf-thm-ae}),
it is standard to verify that
$$
f_{\epsilon}(t,x) \rightarrow \bar f(t,x) \ \ \textrm{a.e. in} \ \{v_* >0\}\ \  \textrm{as}\  \epsilon\rightarrow 0.
$$
Thus due to (\ref{pf-thm-G0}), $f_0=\bar f$ a.e. in $\{v_* >0\}$, i.e.,  $v_*$  satisfies the third equality in (\ref{variationalinequality}).

The proof of Theorem \ref{thm-convergence} is complete.
\end{proof}

\appendix

\section{Proof of Theorem \ref{thm-expansion-boundedness}}

In Appendix A, we present the proof of Theorem \ref{thm-expansion-boundedness}, which is about the expansion and boundedness of
$$
\Omega(t) = \{ x\in\mathbb R^n \ |\ \gamma(t, x)\geq 0 \},
$$
where for convenience, we recall that $\gamma(t,x)$ denotes the solution to the problem  (\ref{one-nonlocalstefan}) with initial data (\ref{initial-onephase}) as follows
\begin{equation}\label{r-one-nonlocalstefan}
\begin{cases}
\displaystyle\gamma_t(t,x)= d \int_{\{\gamma>0\}} \hspace{-0.2cm} k(x-y)\gamma(t,y)dy- d \gamma(t,x)\chi_{\{\gamma>0\}}& t>0, \ x\in\mathbb R^n,\\
\gamma(0,x)= \gamma_0\in C(\bar\Omega_0), \gamma_0|_{\bar\Omega_0}\geq 0, \ \gamma_0|_{\bar\Omega_0} \not\equiv 0& x\in \bar\Omega_0 ,\\
\gamma(0,x)=\gamma_0= -\ell_0  & x\in \mathbb{R}^n\setminus  \bar\Omega_0.
\end{cases}
\end{equation}

\begin{proof}[Proof of Theorem \ref{thm-expansion-boundedness}(i)]
Fix $x\in \mathbb R^n \setminus \bar\Omega_0$. Let $t=s(x)$ denote the moment when $\gamma(s(x), x) = 0$ while $\gamma (t,x)<0$ for $0<t<s(x)$. Set $s(x)=\infty$ if such $s(x)$ does not exist.

By (\ref{one-nonlocalstefan}), when $s(x)<\infty$, one has
\begin{equation*}
\ell_0 = d \int_0^{s(x)}\int_{\mathbb R^n} k(x-y)\gamma^+(\tau, y)dy d \tau = d \int_0^{s(x)}\int_{\Omega(\tau)} k(x-y)
\gamma^+ (\tau, y)dy d \tau.
\end{equation*}
Also thanks to Theorem \ref{thm-wellposedness}, $0\leq \gamma^+ \leq \| \gamma_0|_{\bar\Omega_0}\|_{C(\bar\Omega_0)}$. This yields that
$$
\ell_0  \leq d \int_0^{s(x)}\int_{\Omega(\tau)} k(x-y) \| \gamma_0|_{\bar\Omega_0}\|_{C(\bar\Omega_0)}dy d \tau\leq d s(x)\| \gamma_0|_{\bar\Omega_0}\|_{C(\bar\Omega_0)},
$$
i.e. $s(x)\geq \ell_0 / \left( d \| \gamma_0|_{\bar\Omega_0}\|_{C(\bar\Omega_0)} \right)$.
Hence by choosing $t_0 < \ell_0 / \left( d \| \gamma_0|_{\bar\Omega_0}\|_{C(\bar\Omega_0)} \right)$, one has $\Omega(t) = \Omega(0)$ for $0\leq t\leq t_0$.

The rest follows directly from Proposition \ref{prop-strongmp}.
\end{proof}

In the following, we prove  Theorem \ref{thm-expansion-boundedness}(ii), which is about the uniform boundedness of  $\Omega(t)$.  The proof is lengthy. To control $\Omega(t)$, the main idea is to introduce proper auxiliary $1-$dim problems.

\begin{proof}[Proof of Theorem \ref{thm-expansion-boundedness}(ii)]
To begin with, we introduce the first auxiliary $1-$dim problem
\begin{equation}\label{pf-boundedness-u-M}
\begin{cases}
\displaystyle  \gamma_t(t,x_1)=d  \int_{\mathbb R} k_1 (x_1-y_1)\gamma^+(t,y_1)dy_1- d \gamma^+(t,x_1)   & t>0,\  x_1\in\mathbb{R}, \\
\gamma(0,x_1)=\| \gamma_0|_{\bar\Omega_0}\|_{C(\bar\Omega_0)}    & 0\leq x_1 \leq M,\\
\gamma(0,x_1) =-\ell_0   &  x_1  <0\ \textrm{or}\ x_1>M,
\end{cases}
\end{equation}
where $\displaystyle k_1(x_1) = \int_{\mathbb R^{n-1}} k(x_1, x') dx',\  x' = (x_2,..., x_n)$ and choose  the constant $M$ such that
$$
\bar\Omega_0 \subseteq \{ x\in \mathbb R^n \ |\ 0<x_1<M, \ \textrm{where}\ x=(x_1, ..., x_n)   \}.
$$
Such $M$ exists since $\bar\Omega_0$ is bounded.
Let $\gamma_1 (t,x_1)$ denote the solution to the problem (\ref{pf-boundedness-u-M}).
Notice that $\gamma_1 (t,x_1)$ also
satisfies the $n-$dim problem (\ref{one-nonlocalstefan}) with initial data
$$
\gamma_0(x)=
\begin{cases}
\| \gamma_0|_{\bar\Omega_0}\|_{C(\bar\Omega_0)}    & 0\leq x_1 \leq M,\\
-\ell_0    &  x_1  <0\ \textrm{or}\ x_1>M,\ x=(x_1, ..., x_n).
\end{cases}
$$
Denote
$$
\Sigma_1(t) = \{ x_1\in\mathbb R \ |\  \gamma_1(t, x_1)\geq 0 \}\ \  \textrm{and}\ \  \Sigma_1^{\infty} = \bigcup_{t\geq 0}\ \Sigma_1(t) \subseteq \mathbb R.
$$
By Proposition \ref{prop-cp}, $\gamma_1 (t,x_1) \geq \gamma (t,x)$ in $\mathbb R^n$, where  $x=(x_1, ..., x_n)$,  $\gamma$ denotes the solution to the $n-$dim problem (\ref{r-one-nonlocalstefan}), i.e., (\ref{one-nonlocalstefan}) with initial data (\ref{initial-onephase}).

\smallskip

{\it Obviously,  to prove  Theorem \ref{thm-expansion-boundedness}(ii), it suffices to show that $\Sigma_1^{\infty}$ is bounded,
since the other $n-1$ directions can be handled similarly and thus $\Omega(t)$ will be constrained by a bounded cube.} To verify that $\Sigma_1^{\infty}$ is bounded, for clarity, the rest of the proof is divided into three steps.

\smallskip

{\it Step 1.
We  first show that $|\Sigma_1^{\infty}|$ is bounded.}  Thanks to Proposition \ref{prop-Lp-estimate}, $\gamma_1^+(t, \cdot)\in L^1(\mathbb R)$.  By Proposition \ref{prop-strongmp}, for $0<t<T$, $\Sigma_1(t) \subseteq \Sigma_1(T)$ and thus direct computation yields that
\begin{eqnarray*}
&&\int_{\Sigma_1(T)} \gamma_{1t}(t,x_1) dx_1  \\
  & =& d \int_{\Sigma_1(T)}  \!  \int_{\mathbb R}k_1 (x_1\! -\! y_1)\gamma_1^+(t,y_1) dy_1dx_1 - d \int_{ \Sigma_1(T)} \gamma^+_1(t,x_1) dx_1\\
& \leq& d \int_{\mathbb R}   \gamma_1^+(t,y_1)dy_1 - d \int_{\Sigma_1(T)} \gamma_1^+(t,x_1) dx_1=0.
\end{eqnarray*}
Then, we have
\begin{eqnarray*}
0\leq \int_{\Sigma_1(T)} \gamma_1 (T,x_1) dx_1 \leq\int_{\Sigma_1(T)} \gamma_1 (0,x_1) dx_1= -\ell_0 |\Sigma_1(T)\setminus [0,M]| + \| \gamma_0|_{\bar\Omega_0}\|_{C(\bar\Omega_0)} M,
\end{eqnarray*}
which implies that
\begin{equation*}
|\Sigma_1(T)|   \leq \left(1+ { \| \gamma_0|_{\bar\Omega_0}\|_{C(\bar\Omega_0)}\over \ell_0} \right) M.
\end{equation*}
Since $T$ is arbitrary, one has
\begin{equation*}
| \Sigma_1^{\infty}|   \leq \left(1+ { \| \gamma_0|_{\bar\Omega_0}\|_{C(\bar\Omega_0)}\over \ell_0} \right) M.
\end{equation*}

\smallskip

{\it Step 2. We will show that $\gamma_1(t,x_1)$ decays exponentially as $t$ goes to infinity.}
For this purpose, we introduce the second auxiliary $1-$dim problem with periodic initial data
\begin{equation*}
\begin{cases}
\displaystyle  \gamma_t  =  d  \int_{\mathbb R} k_1 (x_1 - y_1) \gamma^+(t,y_1)dy_1  -  d \gamma^+(t,x_1)   & t>0,\  x_1\in\mathbb{R}, \\
\gamma(0,x_1)=\| \gamma_0|_{\bar\Omega_0}\|_{C(\bar\Omega_0)}   & \kappa (M +L)\leq x_1 \leq \kappa (M +L)+M,\\
\gamma(0,x_1)=-\ell_0    &  \kappa (M +L)+ M < x_1 < (\kappa+1) (M +L),
\end{cases}
\end{equation*}
where $\kappa\in\mathbb Z$, $L>0$ is a constant to be determined later and let $\tilde\gamma_1(t,x_1)$ denote the corresponding solution.
By Proposition \ref{prop-cp},
\begin{equation}\label{pf-periodic-control}
\gamma_1 (t,x_1) \leq \tilde\gamma_1 (t,x_1) \ \ \textrm{for}\  t>0, \, x_1\in\mathbb R.
\end{equation}
Hence, it suffices to verify that $\tilde\gamma_1(t,x_1)$ decays exponentially as $t$ goes to infinity.

First, it is easy to see  that  $\tilde\gamma_1(t,x_1)$ is periodic in $x_1$ with period $M+L$.  Thus this problem can be rewritten as follows
\begin{equation}\label{pf-boundedness-periodic}
\begin{cases}
\displaystyle  \gamma_t  =  d \!  \int_0^{M+L}\! k_* (x_1\! -\! y_1)\gamma^+(t,y_1)dy_1 \! - \! d \gamma^+(t,x_1)   & t>0,\  x_1\in (0,M+L),\\
\gamma(0,x_1)=\| \gamma_0|_{\bar\Omega_0}\|_{C(\bar\Omega_0)}   & 0\leq x_1 \leq M,\\
\gamma(0,x_1)=-\ell_0 <0   &   M < x_1 <  (M +L),
\end{cases}
\end{equation}
where
$$
k_*(x_1) =  \sum_{\kappa\in\mathbb Z} k_1 (x_1 +\kappa (M+L))\ \ \textrm{and}\ \ \int_0^{M+L} k_* (x_1) dx_1 =1.
$$
Denote
$$
\tilde\Sigma_1(t) = \{ x_1\in\mathbb R \ |\  \tilde\gamma_1(t, x_1)\geq 0 \}\ \  \textrm{and}\ \  \tilde\Sigma_1^{\infty} = \bigcup_{t\geq 0}\ \tilde\Sigma_1(t).
$$

Secondly, we need designate the constant $L$.  Fix $T>0$, according to Proposition \ref{prop-strongmp},
$$
[0,M]\subseteq \tilde\Sigma_1(t)\subseteq \tilde\Sigma_1(T),\ \ 0<t<T.
$$
Then, similar to the arguments in {\it Step 1}, by   direct computation, one has for $0<t<T$,
\begin{eqnarray*}
&&\int_{\tilde\Sigma_1(T) \bigcap \,  (0,M+L)} \tilde\gamma_{1t}(t,x_1) dx_1  \\
  & =& d \int_{\tilde\Sigma_1(T) \bigcap \,  (0,M+L)}  \!  \int_0^{M+L}\! k_* (x_1\! -\! y_1) \tilde \gamma_1^+(t,y_1)dy_1dx_1 - d \int_{ \tilde\Sigma_1(T) \bigcap \,  (0,M+L)} \tilde \gamma_1^+(t,x_1) dx_1\\
& \leq& d \int_0^{M+L}  \tilde \gamma_1^+(t,y_1) dy_1 - d \int_{ \tilde\Sigma_1(T) \bigcap \,  (0,M+L)} \tilde \gamma_1^+(t,x_1)  dx_1=0.
\end{eqnarray*}
Thus
\begin{eqnarray*}
0 &\leq & \int_{\tilde\Sigma_1(T) \bigcap \,  (0,M+L)} \tilde\gamma_1 (T,x_1 ) dx_1 \\
& \leq& \int_{\tilde\Sigma_1(T) \bigcap \,  (0,M+L)} \tilde\gamma_1 (0,x_1) dx_1= -\ell_0 | \tilde\Sigma_1(T)\bigcap\, (M,M+L) |  + \| \gamma_0|_{\bar\Omega_0}\|_{C(\bar\Omega_0)} M.
\end{eqnarray*}
This implies that
\begin{equation*}
 | \tilde\Sigma_1^{\infty}\bigcap\, (M,M+L) |  \leq   { \| \gamma_0|_{\bar\Omega_0}\|_{C(\bar\Omega_0)}\over \ell_0} M,
\end{equation*}
since $T$ is arbitrary.
Therefore, we choose $\displaystyle L>{ \| \gamma_0|_{\bar\Omega_0}\|_{C(\bar\Omega_0)}\over \ell_0} M$ to guarantee  that
\begin{equation}\label{a-prop-smaller}
|\, \tilde\Sigma_1^{\infty}\, \bigcap \,  (0,M+L)\, |<M+L.
\end{equation}

Moreover, notice that $\tilde\gamma_1(t, x_1)$ is continuous in $(M,M+L)$. Hence,  thanks to  the strong maximum principle established in  Proposition \ref{prop-strongmp},
one sees that $\tilde\Sigma_1^{\infty}$ is open in $(M,M+L)$.
This observation, together with (\ref{a-prop-smaller}),  implies that there exists a closed  interval $[a,b] \subset (M, M+L)$ satisfying $\displaystyle [a,b]  \bigcap\, \tilde\Sigma_1^{\infty}  = \emptyset$.   If necessary, we could choose $b-a$ smaller such that $k_*(b-a)>0$.
Denote
$$
\tilde\Sigma_D= [0, M]\bigcup (M,a) \bigcup\, (b, M+L).
$$
Obviously, $\tilde\Sigma_1^{\infty} \subseteq \tilde\Sigma_D.$

To estimate $\tilde\gamma_1(t,x_1)$, we first propose  the following  nonlocal eigenvalue problem on the disconnected domain $\tilde\Sigma_D$
\begin{equation}\label{a-EVP}
d \int_{\tilde\Sigma_D} k_*(x_1-y_1)\phi(y_1)dy_1- d\phi(x_1) = \lambda \phi(x_1)\ \ \ \textrm{for} \ x_1\in \tilde\Sigma_D.
\end{equation}
Under the condition that $k_*(b-a)>0$,  the proof of \cite[Theorem 2.6 (i)]{LCW2017} can be  modified to show that the problem (\ref{a-EVP})
admits a principal eigenvalue $\lambda_p$ with the corresponding eigenfunction $\phi_p$ satisfying $\phi_p>0$ in $\tilde{\Sigma}_D$ and also it is easy to see that $\lambda_p< 0$.

Then the third auxiliary problem is introduced as follows
\begin{equation*}
\begin{cases}
\displaystyle v_t(t,x_1)=d \int_{\tilde\Sigma_D} k_*(x_1-y_1)v(t,y_1)dy_1- d v(t,x_1)   & t>0,\  x_1\in\tilde\Sigma_D, \\
v(0,x_1)=\ell \phi_p(x_1)   & x_1\in\tilde\Sigma_D,
\end{cases}
\end{equation*}
where $\ell>0$, and  obviously
$$
v(t,x_1)= \ell e^{\lambda_p t}\phi_p(x_1)
$$
is the corresponding solution.
Choose $\ell$ large enough such that $v(0,x_1) > \| \gamma_0|_{\bar\Omega_0}\|_{C(\bar\Omega_0)}$ in $\tilde\Sigma_D$.
By the comparison principle, it follows that
\begin{equation*}
\tilde\gamma_1 (t,x_1) \leq v (t,x_1) =\ell e^{\lambda_p t}\phi_p(x_1)  \  \ \ \textrm{for}\ \ t>0,\ x_1\in \tilde\Sigma_D.
\end{equation*}
Therefore by (\ref{pf-periodic-control}), the choice of $\tilde\Sigma_D$  and the fact that $\tilde\gamma_1(t,x_1)$ is periodic in $x_1$ with period $M+L$, we have
\begin{equation}\label{pf-thm-compare}
\gamma_1 (t,x_1)\leq  \ell e^{\lambda_p t}\| \phi_p \|_{L^{\infty}(\tilde\Sigma_D)}   \  \ \ \textrm{for}\ \ t>0,\ x_1\in \mathbb R,
\end{equation}
i.e., $\gamma_1(t,x_1)$ decays exponentially at infinity since $\lambda_p< 0$.

\smallskip

{\it Step 3.
We verify that $\Sigma_1^{\infty}$ is bounded, which is the last piece of the proof of Theorem \ref{thm-expansion-boundedness}(ii).}
Suppose that $\Sigma_1^{\infty}$ is unbounded, i.e., there exists a sequence $\{ x_{1i}\}_{i\geq 1}\subseteq \Sigma_1^{\infty}$ and $\{ s_{1i}\}_{i\geq 1}$ with $|x_{1i}|\rightarrow \infty$ as $i\rightarrow \infty$ such that
$$
\ell_0 = d \int_0^{s_{1i}}\int_{\mathbb R} k_1(x_{1i}-y_1) \gamma_1^+(\tau, y_1) dy_1 d \tau,
$$
where $t=s_{1i}$ denote the moment when $\gamma_1(s_{1i}, x_{1i}) = 0$ while $\gamma_1 (t,x_{1i})<0$ for $0<t<s_{1i}$.

To derive a contradiction,  the following property is crucial:
\begin{equation}\label{a-kernel-key}
\lim_{|x_1|\rightarrow \infty} \int_{\Sigma_1^{\infty}}k_1(x_1-y_1)dy_1 =0.
\end{equation}
Since $k_1\in L^1(\mathbb R)$, for any $\epsilon>0$, there exists $R>0$ such that
\begin{equation}\label{a-far}
\int_{|x_1-y_1|\geq R}k_1(x_1-y_1)dy_1 <{\epsilon \over 2}\ \ \textrm{for any}\ x_1\in\mathbb R,
\end{equation}
and there exists $\delta>0$ such that for any measurable set $E \subseteq\mathbb R$ with  $|E|<\delta$, we have
\begin{equation}\label{a-close-small}
\int_{E}k_1(x_1-y_1)dy_1 <{\epsilon \over 2}\ \ \textrm{for any}\ x_1\in\mathbb R.
\end{equation}
Moreover,  according to the property  $| \Sigma_1^{\infty}| <+\infty$ proved in {\it Step 1}, it is standard to show that there exists $X>0$ such that for $|x_1|>X$,
$$
\left | \{|x_1-y_1|< R\}    \bigcap  \{ y_1\in \Sigma_1^{\infty}< \}   \right|<\delta.
$$
This, together with (\ref{a-far}) and (\ref{a-close-small}), yields (\ref{a-kernel-key}) immediately.

Thanks to (\ref{pf-thm-compare}) and (\ref{a-kernel-key}),  it follows that
\begin{eqnarray*}
&& d \int_0^{\infty}\int_{\mathbb R} k_1(x_{1i}-y_1) \gamma_1^+(\tau, y_1) dy_1 d \tau\\
&\leq & d \int_0^{\infty}\int_{\Sigma_1(\tau)} k_1(x_{1i}-y_1)\ell e^{\lambda_p \tau}\| \phi_p\|_{L^{\infty}(\tilde\Sigma_D)}  dy_1 d \tau \\
&\leq &  {d\ell\over - \lambda_p} \| \phi_p\|_{L^{\infty}(\tilde\Sigma_D)}\int_{\Sigma_1^{\infty}} k_1(x_{1i}-y_1) dy_1\rightarrow 0 \ \ \ \textrm{as}\ i\rightarrow \infty.
\end{eqnarray*}
This contradicts to the existence of $s_{1i}$ when $i$ is large enough. Therefore, $\Sigma_1^{\infty}$ is bounded and the desired conclusion is proved.
\end{proof}

\section{Proof of Theorem \ref{thm-continuity}}
Appendix B is devoted to the proof of Theorem \ref{thm-continuity}, which is about the continuous expansion of $\Omega(t)$ under the extra conditions imposed on initial domains and kernel functions.

\begin{proof}[Proof of Theorem \ref{thm-continuity}]
Suppose that the conclusion is not true, i.e., there exists $T>0$ such that  $\Omega(t)$ is connected for $t<T$ while $\Omega(T)$ is disconnected. In other words, $t=T$ is the moment when the jumping phenomena first appears.

Let $\Omega_1(T)\subseteq\Omega(T)$ denote the connected domain which contains $\Omega(t)$ for $t<T$. Choose $y_T\in \Omega(T) \setminus \Omega_1(T)$. $\Omega_1(T)$ exists due to Theorem \ref{thm-expansion-boundedness}(i).
Since $\Omega(0)= \bar\Omega_0$ is convex, there exists a unique $x_T\in \partial\Omega(0)$ such that
$$
|x_T- y_T| = \rm{dist} \{y_T, \Omega(0) \}.
$$
Moreover, there exists $z_T\in\mathbb R^n$, which lies on the line segment $\overline{x_Ty_T}$ and satisfies $z_T \not\in \Omega(T) $. Let $\ell$ denote the line which passes through $(z_T + y_T)/2$ and is perpendicular to the line segment $\overline{x_Ty_T}$. W.l.o.g., assume that $\ell = \{x\in\mathbb R^n \ |\ x_1=0\}$, where $x= (x_1, x_2,..., x_n)$ and $x_{T1}<0$, where $x_T= (x_{T1}, x_{T2},..., x_{Tn})$. Since $\Omega(0)$ is convex, obviously, $\rm{dist} \{\ell, \Omega(0) \}>0$.

For simplicity, denote
$$
\mathbb R^n_- = \{x\in \mathbb R^n\ | \ x_1<0  \},\ \mathbb R^n_+ = \{x\in \mathbb R^n\ | \ x_1>0  \},\ \tilde x = (-x_1, x_2,..., x_n),
$$
and set
$$
w(t,x ) =\gamma (t,x) -\gamma (t, \tilde x),\ x\in \mathbb R^n_-.
$$
Then $y_T=\tilde z_T$  and
\begin{equation}\label{pf-prop-contradiction}
w(T, z_T) = \gamma (t,z_T) -\gamma (t, \tilde z_T)= \gamma (T,z_T) -\gamma (T, y_T)<0.
\end{equation}

Next it is standard to compute that for  $x\in \mathbb R^n_-$,
\begin{eqnarray}
&& w_t(t,x) = \gamma_t (t,x) -\gamma_t (t, \tilde x)\cr
&=&  d\int_{\mathbb R^n} k(x-y)\gamma^+(t,y)dy - d\gamma^+(t,x) -  d\int_{\mathbb R^n} k(\tilde x-y)\gamma^+(t,y)dy  +   d\gamma^+(t,\tilde x)\cr
&=& d\int_{\mathbb R^n_-} k(x-y) \gamma^+(t,y) dy + d\int_{\mathbb R^n_+} k(x-y) \gamma^+(t,y) dy\cr
&& - d\int_{\mathbb R^n_-} k(\tilde x-y) \gamma^+(t,y) dy - d\int_{\mathbb R^n_+} k(\tilde x-y) \gamma^+(t,y) dy - c(t,x) w(t,x)\cr
&=&d\int_{\mathbb R^n_-} k(x-y) \gamma^+(t,y) dy +d \int_{\mathbb R^n_-} k(x-\tilde y) \gamma^+(t,\tilde y) dy\cr
&& - d\int_{\mathbb R^n_-} k(\tilde x-y) \gamma^+(t,y) dy - d\int_{\mathbb R^n_-} k(\tilde x-\tilde y) \gamma^+(t, \tilde y) dy - c(t,x) w(t,x)\cr
&=& \int_{\mathbb R^n_-} \left[ k(x-y ) -k(\tilde x-y) \right] c(t,y) w(t,y)dy - c(t,x) w(t,x),\nonumber
\end{eqnarray}
where
$$
c(t,x) = \frac{d \gamma^+ (t,x) - d   \gamma^+(t,\tilde x)}{\gamma(t,x) - \gamma(t,\tilde x)}.
$$
Note that  $k(x-y ) -k(\tilde x-y)\geq 0$ for $x,y\in \mathbb R^n_-$   since $k(x)$ is decreasing in $|x|$.
Moreover,  for $x\in\ell$, $w(t,x)=0$, and for $x\in \mathbb R^n_-$,
$$
w(0,x ) =\gamma (0,x) -\gamma (0, \tilde x)\geq 0,
$$
since $\Omega(0) \subseteq \mathbb R^n_-$.
Thus by the comparison principle, one has $w(t,x)\geq 0$ for $t>0$, $x\in \mathbb R^n_-$, which contradicts to (\ref{pf-prop-contradiction}).
The proof is complete.
\end{proof}




\bigskip

\noindent\textbf{Data Availability Statement}
\medskip

The authors confirm that this manuscript has no associated data.

\end{document}